\documentclass[lettersize,journal]{IEEEtran}
\usepackage{amsmath,amsfonts}
\usepackage{algorithmic}
\usepackage{algorithm}
\usepackage{array}
\usepackage[caption=false,font=normalsize,labelfont=sf,textfont=sf]{subfig}
\usepackage{textcomp}
\usepackage{stfloats}
\usepackage{url}
\usepackage{verbatim}
\usepackage{graphicx}
\usepackage{cite}
\usepackage{array}
\usepackage{multirow}
\usepackage{tabularx}
\usepackage{graphicx}
\usepackage{booktabs}
\usepackage{makecell}

\newcommand{\cmark}{\checkmark} 
\newcommand{\xmark}{\text{\sffamily X}} 
\usepackage{color}

\begin{document}

\title{An Exact System Optimum Assignment Model for Transit Demand Management}


\author{Xia Zhou, Mark Wallace, Daniel D. Harabor, and Zhenliang Ma,~\IEEEmembership{Member,~IEEE,}
\thanks{This paper was produced by the IEEE Publication Technology Group. They are in Piscataway, NJ.}
\thanks{Manuscript received October 17, 2024; \textit{(Corresponding author: Zhenliang Ma.)}}}

\markboth{Journal of \LaTeX\ Class Files,~Vol.~14, No.~8, October~2024}%
{Shell \MakeLowercase{\textit{et al.}}: A Sample Article Using IEEEtran.cls for IEEE Journals}


\maketitle

\begin{abstract}
Mass transit systems are experiencing increasing congestion in many cities. The schedule-based transit assignment problem (STAP) involves a joint choice model for departure times and routes, defining a space-time path in which passengers decide when to depart and which route to take. User equilibrium (UE) models for the STAP indicates the current congestion cost, while a system optimum (SO) models can provide insights for congestion relief directions. However, current STAP methods rely on approximate SO (Approx. SO) models, which underestimate the potential for congestion reduction in the system. The few studies in STAP that compute exact SO solutions ignore realistic constraints such as hard capacity, multi-line networks, or spatial-temporal competing demand flows. The paper proposes an exact SO method for the STAP that overcomes these limitations. We apply our approach to a case study involving part of the Hong Kong Mass Transit Railway network, which includes 5 lines, 12 interacting origin-destination pairs and 52,717 passengers. Computing an Approx. SO solution for this system indicates a modest potential for congestion reduction measures, with a cost reduction of 17.39\% from the UE solution.  Our exact SO solution is 36.35\% lower than the UE solution, which is more than double the potential for congestion reduction. We then show how the exact SO solution can be used to identify opportunities for congestion reduction: (i) which origin-destination pairs have the most potential to reduce congestion; (ii) how many passengers can be reasonably shifted; (iii) future system potential with increasing demand and expanding network capacity. 

\end{abstract}


\section{Introduction}
\IEEEPARstart{M}{ass} transit systems are facing increasing crowding due to limited service capacity and rising travel demand \cite{ref1}. One way to alleviate system crowding is to effectively manage demand, maximizing the utilization of existing capacity \cite{ref2}. 
An efficient and reliable demand assignment method is crucial for policymakers to evaluate the potential congestion relief of the current system \cite{ref3}. 
Insights into the features of the schedule-based transit assignment problem (STAP) are crucial for understanding the transit system.
The STAP assigns the demand flows over time and/or space under a certain objective, such as the user equilibrium (UE) or the system optimum (SO). The UE assignment is characterized by passengers' free and selfish behavior, in which each individual seeks to minimize their own travel costs. This UE approach is useful for for approximating a system state without external interventions. Unlike the UE assignment, the SO assignment assumes that passengers cooperate to minimize the overall system cost. This SO approach is useful for representing the ideal system with the best performance. A deeper understanding of the UE and SO modeling in the transit system is crucial to guide the policy making process, as it provides information on which passengers contribute the most to the system  congestion and the potential benefits for policy interventions. 
\begin{figure}[!t]
\centering
\includegraphics[width=3.2 in]{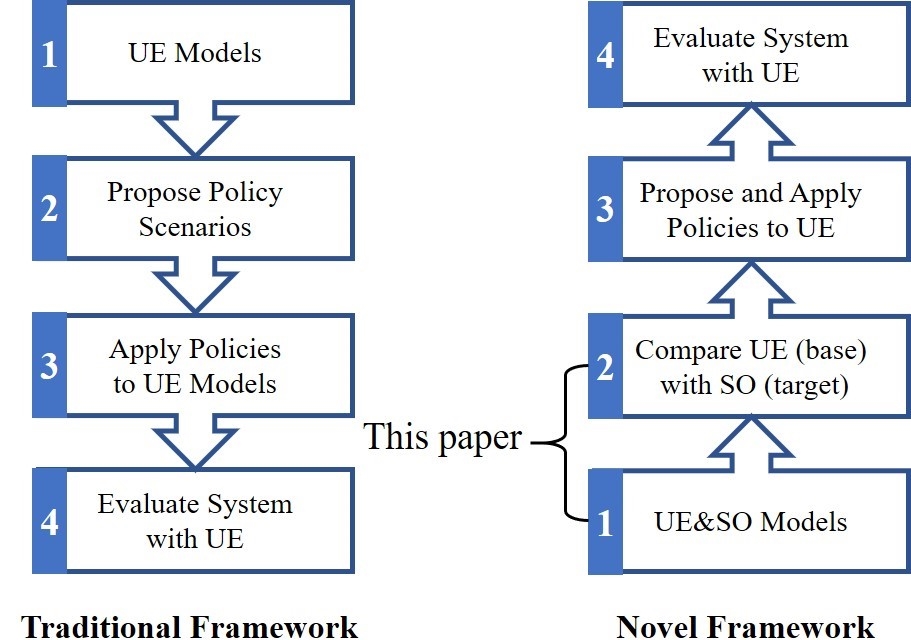}
\caption{Frameworks for evaluating promotion policies in the transit system (adapted from \cite{ref56})}
\label{Fig16}
\end{figure}

As shown in Fig. \ref{Fig16}, in the traditional evaluation method, a congestion relief policy or promotion scheme is firstly proposed, and then its impact is evaluated using mathematical tools, typically UE transit assignment models (for example, \cite{ref54, ref55}). In the traditional framework, transit operators and planners cannot determine whether a proposed policy will achieve optimal outcomes, as the upper bound potential of the system (the SO system) is absent. The SO solution offers insights into cost descent directions and shows the bound for potential cost reductions. Instead of using the traditional framework, \cite{ref56} proposed a 'reverse engineering' framework, which compares the UE system (base) with the SO system (target) to deduce effective policies. With the upper bound of the system available, transit operators and planners can develop informed policies to alleviate congestion. Unfortunately, computing SO transit assignment solutions is extremely challenging.


\begin{figure*}[ht]
  \centering
  \begin{minipage}[b]{0.46\textwidth}
    \centering
    \includegraphics[width=\textwidth]{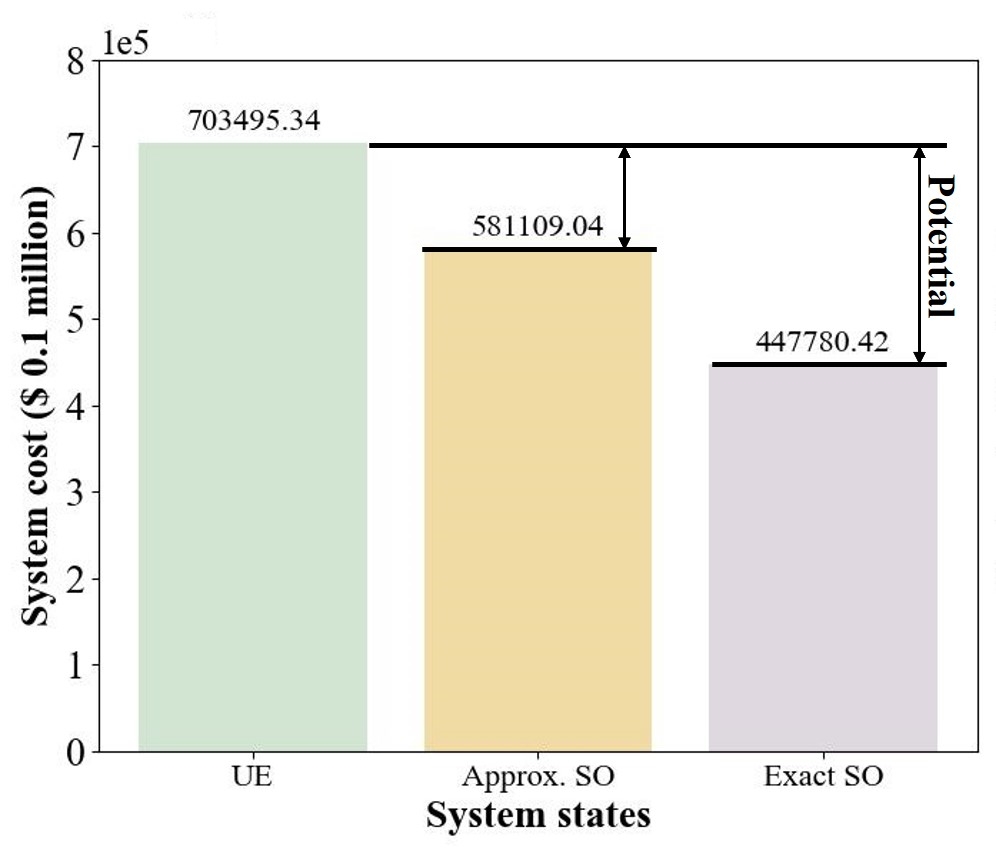} 
    \par\noindent
    \textbf{(a)}
    \label{fig2:minipage1}
  \end{minipage}
  \hfill
  \begin{minipage}[b]{0.46\textwidth}
    \centering
    \includegraphics[width=\textwidth]{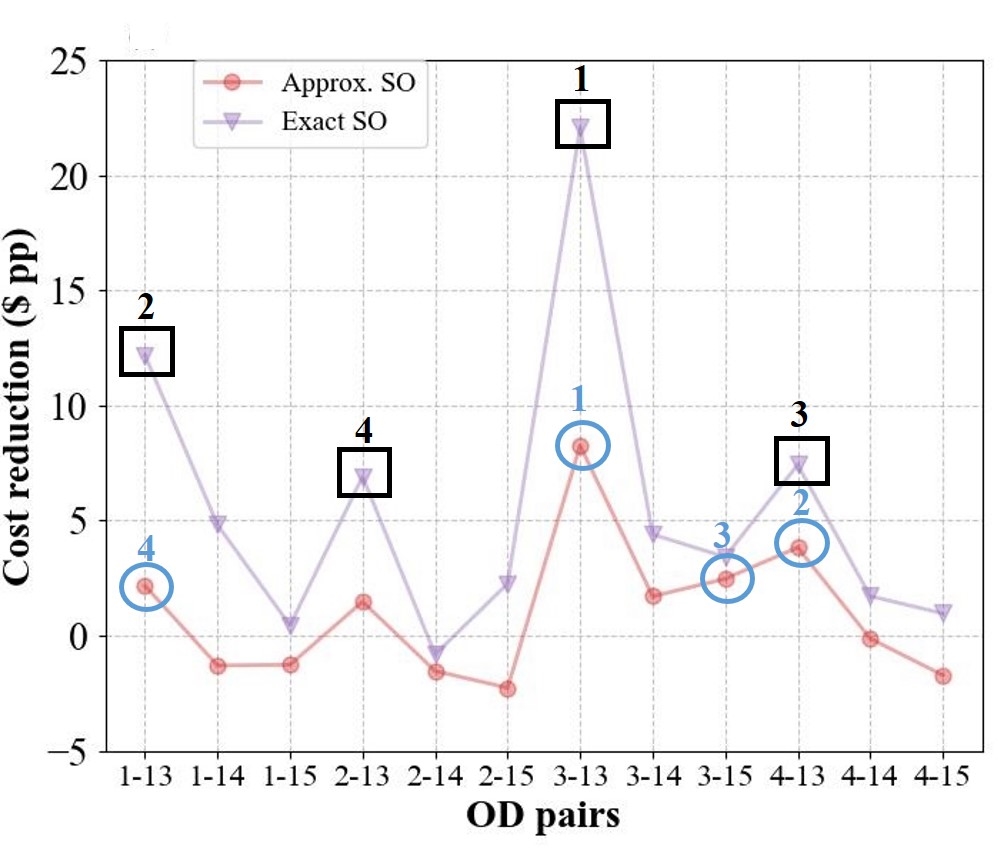} 
    \par\noindent
    \textbf{(b)}
    \label{fig2:main}
  \end{minipage}
  \caption{Comparison between Approx. SO and exact SO: (a) Measured system potential for promotions (current: UE, target: SO); (b) Potential cost reduction contribution per passenger: targeting the top 4 groups for potential shifts (a higher cost reduction indicates greater potential for congestion relief)}
  \label{fig2}
\end{figure*}


Approaches for solving the system-optimal STAP problem can be broadly categorized into approximate SO (Approx. SO) methods (for example, heuristic, numerical, simulation-based methods) and exact SO solution methods (analytical methods). In general, Approx. SO methods have two main drawbacks: (a) They provide local optimum solutions instead of global ones, with solution quality/cost highly dependent on the algorithm used and available computational resources; (b) They may provide misleading information (for example, targeting the wrong passengers for shifting, underestimating the potential improvement of the system) for guiding congestion relief. 
To the best of the authors' knowledge, existing studies on exact SO STAP are extremely limited. No study was found to consider the simultaneous inclusion of hard train capacity, multi-line networks, and spatial-temporal flows (encompassing both departure times and route choices). 
Figure \ref{fig2} illustrates the guidance provided by both Approx. and exact SO models, including system potential evaluations for congestion relief and the identification of the top potential passenger groups for shifting. Fig. \ref{fig2}(a) shows that the Approx. SO underestimates the potential improvements of the UE system. Fig. \ref{fig2}(b) shows that the Approx. SO may target the wrong passenger groups for shifting. Detailed results can be found in the case study section.

This study proposes an exact SO model for the schedule-based transit assignment problem to provide a global optimal SO solution. The proposed approach is evaluated on a subset of the Hong Kong Mass Transit Railway (MTR) network, using real-world demand data. In summary, this paper contributes to the state of the art in the system-optimal STAP problem with the following advancements:
\begin{itemize}
    \item Propose the exact SO model of the STAP, considering the hard train capacity in multi-line networks. Compared to the results of Approx. SO, the exact SO provides more accurate and reliable information for evaluating the policy potential.
    \item Using the UE and exact SO model developed in this paper, we derive managerial insights and policy suggestions for the Hong Kong MTR network under different scenarios. The scenarios include (a) varying available passenger choices, and (b) varying demand levels and capacity levels, corresponding to the future development of the city.
\end{itemize}

The remainder of the paper is organized as follows. Section \ref{sec:literature-review} reviews relevant literature on the STAP. Section \ref{sec:problem-formulation} introduces the problem description of the STAP. Section \ref{sec:UE-approx-SO} formulates 
the UE and SO models, which are evaluated and solved in Section \ref{sec:solution-procedure}. Section \ref{sec:case-study} validates the proposed exact SO model using an extracted Hong Kong MTR network. The last Section \ref{sec:conclusion} summarizes the main findings and future research.

\section{Literature review} \label{sec:literature-review}
Exact SO models for the traffic assignment problem (TAP) have been widely studied in car traffic (i.e., \cite{ref4, ref5, ref6, ref27, ref29, ref30, ref32}). 
Analytical models were proposed in \cite{ref5} for the system-optimal traffic assignment problem in car traffic, aiming to solve for the time-dependent link and path flows within a network. 
Unlike the open car traffic system, the closed system of the transit system is constrained by predefined routes and schedules \cite{ref7}. 
These bring additional modeling challenges including, for example, transfer times between service lines, and times waiting on a platform for the next train when the first train is too full to board due to constrained vehicle capacity.
Importantly, the discrete, nonlinear nature of system dynamics and travel costs makes the STAP more complex compared to the TAP in the continuous-flow car traffic system. Therefore, TAP models in car traffic are not easily applied or extended to the transit system. In this paper, we propose an exact SO model for the STAP considering simultaneous departure time choice, route choice, hard capacity, and multi-line networks.

\begin{table*}[ht]
\caption{System optimum assignment in the transit system- Departure time choice and Route choice}
\label{tab:table1}
\centering
\begin{tabularx}{0.86\textwidth}{p{1.2cm} p{1.0cm} p{1.0cm} p{1.6cm} p{3.0 cm} p{1.2cm} p{1.4cm}p{1.4cm}}
\toprule
\multirow{2}{*}{\textbf{Ref.}} & \multicolumn{2}{c}{\centering \textbf{Decision Variables}} & \centering \textbf{Network} & \centering \textbf{Constraint} & \multicolumn{3}{c}{\textbf{Objectives}}\\
\cmidrule(lr){2-3} \cmidrule(lr){6-8}
 & \centering DT & \centering RC & \centering (Multi-line) & \centering (Hard Train Capacity) & UE & Approx. SO &  \textbf{Exact SO} \\
\midrule
\multicolumn{8}{c}{\textit{Papers considering Approx. SO (DT or RC)}} \\
\midrule
\cite{ref8} & \centering \xmark & \centering \cmark  &  \centering \xmark & \centering \cmark & \xmark & \cmark & \xmark\\
\cite{ref18} & \centering \cmark & \centering \xmark & \centering \xmark & \centering \cmark & \cmark & \cmark & \xmark\\
\cite{ref19} & \centering \cmark & \centering \xmark & \centering \xmark & \centering \xmark  & \cmark & \cmark & \xmark\\

\cite{ref22} & \centering \xmark & \centering \cmark  & \centering \cmark & \centering \cmark & \cmark & \cmark & \xmark\\
\midrule
\multicolumn{8}{c}{\textit{Papers considering Approx. SO (DT \& RC)}} \\
\midrule
\cite{ref25} & \centering \cmark & \centering \cmark  & \centering \xmark & \centering \cmark & \xmark & \cmark & \xmark\\
\cite{ref26} & \centering \cmark & \centering \cmark  & \centering \cmark & \centering \xmark & \cmark & \cmark & \xmark\\
\midrule
\multicolumn{8}{c}{\textit{Paper(s) considering Exact SO}} \\
\midrule
\cite{ref9} & \centering \cmark & \centering \xmark  & \centering \xmark & \centering \xmark & \cmark & \xmark & \cmark\\
\midrule
\multicolumn{8}{c}{\textit{This paper}} \\
\midrule
This paper & \centering \cmark & \centering \cmark  & \centering \cmark & \centering \cmark & \cmark & \cmark & \cmark\\
\multicolumn{8}{l}{} \\
\multicolumn{8}{l}{\textbf{DT:} Departure Time (Choice); \textbf{RC:} Route Choice; \textbf{UE:} User Equilibrium; \textbf{Approx. SO:} Approximate System Optimum; } \\
\bottomrule
\end{tabularx}
\end{table*}
In the literature, models for the transit assignment problem are typically categorized as either frequency-based 
or schedule-based 
transit assignment models \cite{ref51}. The frequency-based transit assignment (FTA) model aims to determine the average loads on transit lines and to identify potential macroscopic congestion \cite{ref52}. In FTA models, the time dimension is absent, and all model characteristics are averaged over the entire time period. Therefore, they are unable to accurately evaluate the transit network state under changing conditions, such as temporal variations in passenger flows. The STAP model aims to determine the loads on individual service runs and their spatial-temporal flow distributions \cite{ref53}. The growing focus on managing spatial-temporal flows in interconnected transit networks by optimizing transfers and minimizing passenger delays clearly necessitates the use of STAP models. In this paper, we focus on STAP models.

The routing strategy (how the route input is generated) is the core element of the the STAP.
Approaches for routing strategies include the following \cite{ref52}: (1) Agent-based simulation; (2) Spatial-temporal network (time-expanded network).
The agent-based method is highly flexible and better suited for complex, realistic, large-scale transit assignment problems, such as those involving travel strategies and seat availability \cite{ref51}, or demand uncertainty \cite{ref8}. However, the agent-based simulation comes with the drawback of producing unstable and approximate results. The spatial-temporal expanded  method is commonly used in mathematical
models, directly yielding the expected values of output variables such as loads and travel costs.  See, for example, \cite{ref24}. 
This paper focuses on exact SO models for the STAP. Therefore, we use the spatial-temporal expansion method. 



In the transit system, as passengers' travel cost decisions are interconnected and there are hard capacity constraints, travel cost functions are unpredictable \cite{ref8}.
It is a challenging task to formulate the exact SO model for the STAP, considering simultaneous hard train capacity, a multi-line network with transfer behaviors and interacting OD flows.

Given these features, solving the STAP remains a challenge. Table \ref{tab:table1} shows that most transit studies either neglect capacity (i.e., \cite{ref19}) or use soft  capacity constraints, 
assuming few passengers are left on the platform when a train departs, as they squeeze in regardless of the train's capacity.  Furthermore, most studies propose Approx. SO models for the STAP considering departure time choice \cite{ref18}, route choice \cite{ref22}, simultaneous departure time choice and route choice \cite{ref25}, multi-line networks \cite{ref26}, while only one paper \cite{ref9} proposed an exact SO model for the STAP. That study formulated an exact SO approach to manage peak-hour congestion on a single rail transit line. However, it does not account for hard train capacity and route choices. To the best of the authors' knowledge, no study has proposed an exact SO model for the STAP that simultaneously considers departure time choice, route choice, hard capacity, and multi-line networks.

SO solutions that simultaneously consider these features are more realistic compared to those in existing studies. More realistic solutions will offer more valuable insights for guiding congestion relief. However, incorporating too many detailed features, especially less important ones, may make the problem hard to solve. Therefore, this paper focuses on the main (realistic) features: departure time choice, route choice, hard capacity, and multi-line networks.

\section{Problem Formulation}\label{sec:problem-formulation}
.
\subsection{The generic transit assignment model}
The model comprises a transit network, a public transport timetable, a set of passenger demands, a set of journey constraints and a set of costs. In our case, the transit network comprises a set of train lines, with stations on each line, some of which are transfer stations between the lines.

A sequence of trains run on each line according to a timetable which gives a fixed departure time for each train from each station.  
We make the following assumptions as reported in \cite{ref3, ref8, ref34}:
\begin{itemize}
\item In-vehicle travel times are constant;
\item Schedules of trains are fixed, and train dwell times are included in schedules;
\item Train capacity is fixed, meaning passengers must wait for the next train when the current one is full.
\end{itemize}

There is a set of origin-destination (OD) passenger demands, where each origin and each destination is a station.  All passengers have a target arrival time, such as 9:00 AM.

The cost to a passenger is a sum of weighted costs of journey time components.  The costs comprise waiting time, in-vehicle time, early arrival time, and late arrival time \cite{ref9, ref33}. 

This generic transit assignment model covers three model variants, including (1) User Equilibrium (UE) STAP with users choosing routes and departure times, (2) Approx. SO STAP using simulation for network loading and (3) Exact SO STAP using analytic models for network loading. 

The UE approach uses simulation for the network loading to reflect how each passenger optimises their own travel cost by choosing his/her route and departure time.  
 
 
 The Approx. SO model also uses simulation for the network loading, but in this case passenger choices are adapted to optimise the total system cost. The simulation used by the Approx SO models passengers boarding (or being denied boarding) in the same way as the UE simulation.  The difference is in the passengers' choice of departure time and route.

 
 
The Exact SO model optimizes the total system cost by determining each passenger's route and departure time, assuming that in the SO, every passenger catches the next train on their route departing from their current station (`no denied boarding' holds in transit systems with sufficient trains provided). The model is sufficiently simplified to be solved to optimality without the need for simulation. 
 

\subsection{Model Settings and Notations}
For both UE and SO travel patterns, the model must include the makeup of passengers who travel on each train and each leg of its route.  In the UE, Approx. SO and Exact SO, the passengers are grouped by OD and route, and the number of passengers in each group on each leg of each train are modelled.  

Passengers' decisions on  which departure train ($t$) and which route ($r$) at their origin ($o$) are denoted as ($t, r$). All notations used in this paper are summarized in Tables \ref{tab:notations-parameters} and \ref{tab:notation-variables}.
\begin{table}[ht]
    \centering
    \begin{minipage}[b]{0.46\textwidth}
        \centering
        \caption{Notations - Model Parameters and Variables}
        \label{tab:notations-parameters}
        \begin{tabularx}{\columnwidth}{cX X}
            \toprule
            \multicolumn{2}{l}{\textbf{Network Parameters}}\\
            \midrule
            $L$ & the set of lines \\
            $S$ & the set of all stations\\
            $S_l$ & index the stations on line $l$ ($S_l = 1..sct_l$ where $sct_l$ is the number of stations on line $l \in L$) \\
            $s^m_l$  & the $m^{th}$ station on line $l$\\
            $e^m_l$  & the $m^{th}$ leg (`e' is for `edge') on $l$ , between stations $s^m_l$ and $s^{m+1}_l$\\ 
            $s$ & is the transfer station if $s = s^{m_1}_{l_1} = s^{m_2}_{l_2}$ where $l_1 \neq l_2$; Otherwise it is the normal station\\
            $ne(s,l)$ & the edge following station $s$ on line $l$ (exiting edge)\\
            $pe(s,l)$ & the edge preceding $s$ on line $l$ (entering edge) \\
            $l^t$ & the line of train $t$ on\\
            \midrule            
            \multicolumn{2}{l}{\textbf{Timetable Parameters}}\\
            \midrule
            $T$ & the set of trains \\
            $cap(t)$ & the capacity of train $t \in T$\\
            $T_l$ & index the trains on line $l$, $T_l = 1..tct_l$ where $tct_l$ is the number of trains on line $l \in L$\\
            $t^n_l$ & the $n^{th}$ train on line $l$\\
            $tt^n_l$ & the time at which train $t^n_l$ departs from the first station $s^1_l$\\
            $dur(e)$ & time for a train to traverse leg $e$\\
            $dep(t^n_l,s^m_l)$ & the departure time for train $t^n_l$ leaves station $s^m_l$ ($dep(t^n_l,s^m_l) = tt^n_l + \sum_{i<m} dur(e^i_l)$)\\
            $arr(t,s)$ & the arrival time for train $t$ at station $s$ (we assume the arrival time is same as the departure time: $arr(t,s) = dep(t,s)$)\\
            $con(s, t_1,t_2)$ & connection between two trains (Boolean algebra): a connection holds (is true) if $s$ is a transfer station for line $l_1$ and $l_2$ and if $t_2$ is the first train to depart on line $l_2$ after the arrival of $t_1$ on $l_1$.\\
            \midrule            
            \multicolumn{2}{l}{\textbf{Passenger Demand Parameters}}\\
            \midrule
            $K$ &  the set of OD pairs with nonempty demand \\   
            $o^k$ & the origin of OD $k \in K$\\
            $O$ & the set of origins $\{o^k : k \in K\}$ \\
            $d^k$ & the destination of OD $k$\\
            $D$ & the set of destinations $\{d^k : k \in K\}$\\
            $q^k$ & the number of passengers for OD $k$ \\
            $R^k$ & the set of possible routes for OD $k$\\
            $l_{o^k}^r$ & the line from $o^k$ used by route $r \in R^k$\\  
            $ch(r,s,l_1,l_2)$ & connection between two lines (Boolean algebra): is true if on route $r$, passengers transfer at station $s$ from line $l_1$ to $l_2$\\
            $l_{d^k}^r$ & the line to $d^k$ used by route $r$\\
            \midrule            
            \multicolumn{2}{l}{\textbf{Cost Parameters}}\\ 
            \midrule
            $\alpha$ & the cost weighting for in-vehicle travel time (i.e., the cost for a passenger traversing edge $e$ is $\alpha \times dur(e)$)\\
            $\beta$ & the cost weighting for time waiting on the platform (i.e., the cost for a passenger making connection $con(s, t_1, t_2)$ is $\beta \times (dep(t_2,s) - arr(t_1,s))$\\
            $\gamma$ & the cost weighting for arriving too early at the destination \\
            $\mu$ & the cost weighting for arriving too late at the destination \\
            \bottomrule
        \end{tabularx}
    \end{minipage}
\end{table}

\begin{table}[ht]
    \begin{minipage}[b]{0.46\textwidth}
        \centering
        \caption{Notations - Model Variables}
        \label{tab:notation-variables}
        \begin{tabularx}{\columnwidth}{cX X}
            \toprule          
            \multicolumn{2}{l}{\textbf{Model Flow Variables}}\\ 
            \midrule
            $q^k(t, r)$ & the number of passengers whose OD is $k$, that arrive at their origin station to meet train $t$ setting out on route $r \in R^k$\\
            $f^k_e(t, r)$ & the number of passengers on route $r$ that are on train $t$ and leg $e$ \\
            $g^k_s(t, r)$ & the number of passengers on route $r$ at station $s$ awaiting train $t$\\
            $db^k_s(t, r)$ & the number of passengers on route $r$ at station $s$ waiting for train $t$ because they could not board the previous train\\
            $fk^k_e(t)$ & the total number of passengers whose OD is $k$ that are on train $t$ and leg $e$\\
            $dbk_s^k(t)$ & the total number of passengers whose OD is $k$ at station $s$ waiting for train $t$ because they could not board the previous train\\
            $I^k(t, r)$ & index passengers whose OD pair is $k$ on route $r$ meeting train $t$ at origin $o^k$\\ 
            $\mathbf{\textit{Q}}$ & passenger demand distribution over trains and routes at origins: $\mathbf{\textit{Q}} = \{q^k(t, r): k\in K, r \in R^k, t \in T_l$, where $l = l^r_{o^k}\}$\\
            $\mathbf{\textit{F}}$ & passenger flow distribution over legs connecting to destinations:
            $\mathbf{\textit{F}} = \{f^k_{pe(d^k, l)}(t, r): k \in K, r \in R^k, t \in T_l$, where $l = l^r_{d^k}\}$\\
            \midrule
            \multicolumn{2}{l}{\textbf{Model Cost Variables}}\\ 
            \midrule
            $ivcst1_e^k(t)$ & the in-vehicle cost of the passengers whose OD is $k$ on train $t$ and leg $e$\\
            $ivcst(k)$ & the total in-vehicle cost of the passengers whose OD is $k$\\
            $chcst1^k_s(t_1,t_2)$ & the waiting cost of the passengers whose OD is $k$ at station $s$ changing from train $t_1$ to train $t_2$\\
            $chcst(k)$ & the total waiting (changing between two adjacent trains) cost of passengers whose OD is $k$\\
            $dbcst(k)$ & the total waiting (denied boarding) cost of passengers whose OD is $k$\\
            $elcst1^k(l,n)$ & the late/early penalty cost of passengers whose OD is $k$ arriving at their destination $d^k$ on the $n^{th}$ train of line $l$\\
            $elcst(k)$ & the total late/early penalty cost of passengers whose OD is $k$\\
            $cost(k)$ & the total cost of passengers whose OD is $k$\\
            $c^k_i(t,r)$ & the travel cost of the $i^{th}$ passenger in $I^k(t, r)$\\
            $avc^k(t,r)$ & the average cost of option ($t, r$) chosen by passengers in $I^k(t, r)$ at origin $o^k$\\
            $avc^k_0(t,r)$ & the free flow cost of option ($t, r$)\\
            $avc^k_{*}$ & the average cost of the best option within OD $k$\\
            \bottomrule
        \end{tabularx}
    \end{minipage}

\end{table}
The model also includes the number of passengers, in each group, waiting at the stations on each line.  When the next train arrives on the line, the model determines which passengers board the train, and which, if any, are left on the platform because the train is full. 

In principle, passengers who arrive first on the platform should be the first to board the train. For passengers arriving at the same time from different locations, models typically capture a mix of these passengers and determine who can board the next train.
This feature of the model is supported by simulations in UE and Approx. SO models.

We make a key novel observation here.  
If all the passengers on the platform can board the train, there is no need for simulation to support this aspect of the model.
Based on this observation we introduce 
the model for Exact SO, in which all passengers can board the next train.
\ \\
\subsection{Model constraints}
\label{sec:constraint}
For readability, the presentation is somewhat informal. The complete formal constraint specification is given in the appendix (see Appendix A).

The 'train capacity' constraint limits the total number of passengers on each train $t$ on each leg $e$:\\
\begin{equation}
    \sum{k \in K,r \in R^k}: \  f^k_e(t,r) \leq cap(t) 
    \label{eq:1}
\end{equation}

The next group of constraints model how many passengers are (on the platform) waiting for a train.

At the origin of an OD, the passengers waiting on the platform are those selecting that train to start, plus those who could not catch the previous train:

If $s=o^k$ then
\begin{equation}
   g^k_s(t, r) = q^k(t,r) + db^k_s(t, r)
   \label{eq:2}
\end{equation}
Otherwise, at a transfer location for the route $r$, the passengers waiting on the platform are those that changed from a connecting train on the previous line on route $r$, plus those who missed the previous train:

If $s \neq o^k$ and 
train $t = t^n_{l_2}$ and if 
$ch(s,r,l_1,l_2)$ (is true) and $con(s,t_1,t)$ (is true) then:\\
\begin{equation}
g^k_s(t, r) = f^k_{pe(s,l_1)}(t_1,r) + db^k_s(t,r)
\label{eq:3}
\end{equation}

Note that unless the route includes a cycle, no passengers on route $r$ can already have been on the train when it arrived at $s$.
Consequently, $f^k_{pe(s,l_1)}(t_1, r)$ is entirely made up of passengers from $g^k_s(t, r)$ who were waiting at the station.

If $s$ is not the origin nor a transfer station  for route $r$, the number of waiting passengers is $0$:

If $s \neq o^k$ and not $ch(s,r,l_1,l_2)$ then
\begin{equation}
g^k_s(t, r)=0
\label{eq:4}
\end{equation}

Also, if $t$ is the first train on the line, then no passengers are waiting who missed the previous train.
Thus, for all ODs, routes, lines and stations:
\begin{equation}
db^k_s(t^1_l, r) =0
\label{eq:5}
\end{equation}

Similarly, everybody must board the last train, so:
\begin{equation}
db^k_s(t^{tct_l+1}_l, r) = 0
\label{eq:6}
\end{equation}

The next group of constraint models the passengers on train $t$ traversing leg $e$.

Firstly, if the route $r$ does not have a transfer at $s$, and $s$ is not the final station on the line, or origin or destination, then all passengers continue on the train:
\begin{equation}
f^k_{pe(s,l)}(t, r) = f^k_{ne(s,l)}(t, r)
\label{eq:7}
\end{equation}
unless $s=s^1_l$ or $s=o^k$ or $s=d^k$ or $ch(s,r,l,l_2)$ where $l_2 \neq l$. Otherwise, the number of passengers getting on train $t$ (on line $l^t$) is no more than the number waiting on the platform: 
\begin{equation}
f^k_{ne(s, l^t)}(t,r) \leq g^k_s(t, r)
\label{eq:chboard}
\end{equation}

However the actual number that board the train - and the number who fail to board it - is predicted by a simulation rather than a mathematical constraint.

The next group of constraints imposes that all passengers on a given route eventually depart from their origin and arrive at their destination.

Departing passengers:

If $l = l_{o^k}^r$ then
\begin{equation}
q^k = \sum_{t \in T_l} \sum_{r \in R^k} q^k(t, r) = \sum_{t \in T_l} \sum_{r \in R^k} f^k_{ne(o^k,l)}(t, r)
\label{eq:9}
\end{equation}

Arriving passengers:

If $e$ is the leg before $d^k$, $l_1 = l_{o^k}^r$, and $l_2 = l_{d^k}^r$,  then:
\begin{equation}
q^k = \sum_{t_1 \in T_{l_1}} \sum_{r \in R^k} q^k(t_1, r) = \sum_{t_2 \in T_{l_2}} \sum_{r \in R^k} f^k_e(t_2, r)
\label{eq:10}
\end{equation}

The final group of constraints specifies the cost.  We divide the cost into three parts: in-vehicle cost; waiting cost; and early/late arrival cost.

For the following constraints we use:\\
\begin{equation}
fk^k_e(t) = \sum_{r \in R^k} f^k_e(t,r)
\label{eq:100}
\end{equation}
and
\begin{equation}
dbk_s^k(t) = \sum_{r \in R^k} db^k_s(t, r)
\end{equation}
for total passengers with OD $k$ on that are on train $t$ and leg $e$, and total denied boarding passengers with OD $k$ at station $s$ waiting for train $t$ respectively.

The in-vehicle cost is the sum of the passenger loads on each leg.  For clarity, for the cost of the load on each train and leg, we write:
\begin{equation}
ivcst1_e^k(t) = \alpha \times fk^k_e(t) \times dur(e)
\label{eq:121}
\end{equation}

The in-vehicle cost is:\\
\begin{equation}
ivcst(k) = \sum_{l \in L,n \in T_l,m \in S_l} : ivcst1^k_{e^m_l}(t^n_l)
\label{eq:11}
\end{equation}

The waiting cost is the sum of the cost due to changing trains.
The cost $chcst$ of a single transfer is specified as follows.
If  
train $t_2 = t^n_{l_2}$, 
$ch(s,r,l_1,l_2)$ and $con(s,t_1,t_2)$ then:
\begin{flalign}
chcst1^k_s(t_1,t_2) = \\
\notag \sum_{r \in R^k} : \beta &\times (dep(t_2,s)-arr(t_1,s))\\ 
\notag &\times f^k_{pe(s,l_1)}(t_1, r) 
\label{eq:122}
\end{flalign}

The cost for all transfers where $ch(s,r,l_1,l_2)$ holds is:\\
\begin{flalign}
chcst(k) = &\sum_{s \in S, t_1 \in T_{l_1}, t_2 \in T_{l_2}} : chcst1^k_s(t_1,t_2)
\end{flalign}

The cost for passengers denied boarding is 
\begin{flalign}
dbcst(k) = \sum_{l \in L, n \in S_l} & \beta \times (tt^n_l-tt^{n-1}_l) 
\\ \notag & \times dbk^k_s(t^n_l) 
\end{flalign}

The late/early penalty costs are, for a train arriving at a destination before or after 9:00 AM (writing $540$ for 9:00 AM).
If $t = t^n_l$ and $d^k = s^m_l$ for some $m \in S_l$ then:
\begin{flalign}
\label{eq:123}
elcst1^k(l,n) = & fk^k_{pe(d^k, l)}(t) \times\\
\notag & (\gamma \times \max(0, 540-arr(t,d^k)) +\\ 
\notag & \mu \times \max(0, arr(t,d^k)-540))
\end{flalign}

The total late/early penalty is:\\
\begin{equation}
elcst(k) = \sum_{l \in L, n \in T_l} : elcst1^k(l,n)
\end{equation}

Finally:
\begin{equation}
    cost(k) = ivcst(k)+chcst(k)+dbcst(k)+elcst(k)
    \label{eq:costr}
\end{equation}

\subsection{Simulated travel costs for UE and Approx. SO}

There are components of the model which require simulation when predicting the UE and the SO.  These are associated with Eq. (\ref{eq:chboard}).
In each case, simulation is required to determine which passengers catch the train. The values of the other variables are also computed by the simulator using the Eqs. (\ref{eq:1})-(\ref{eq:10}).  

The simulation model (network loading) is based on the event-based network performance model structure in \cite{ref48} and \cite{ref49}. It works as follows: the events in the system include the arrival and departure of train runs at stations.  Events are sequenced in time.  At each train departure, the passengers on the train comprise: (1) the passengers already on the train before this station, and (2) a fraction of the passengers waiting at the station platform for this train: the fraction is 1 if there is enough room on the train, otherwise it is the fraction that fills the train to capacity. 
At each train arrival, the passengers leaving the train are: (1) those who have arrived at their destination, and (2) those who are changing to the next train on their routes.  These passengers proceed to the relevant platform and wait for the next train. The input to the simulation is $\mathbf{\textit{Q}} = \{q^k(t, r): k\in K, r \in R^k, t \in T_l$, where $l = l^r_{o^k}\}$.  Its output is $\mathbf{\textit{F}} = \{f^k_{pe(d^k, l)}(t, r): k \in K, r \in R^k, t \in T_l$, where $l = l^r_{d^k}\}$. The cost is then calculated using Eqs. (\ref{eq:100})-(\ref{eq:costr}).

\section{User Equilibrium and System Optimum Models}\label{sec:UE-approx-SO}
\subsection{User Equilibrium}
The formulation of UE assignment usually follows the extension of Wardrop's first principle \cite{ref7}: the perceived travel cost for users across all available options for an OD pair is equal and minimized. Conversely, the perceived travel cost for unused options is higher than or equal to that of the utilized options.

Such UE conditions were formulated as a gap-minimizing mathematical model in car traffic \cite{ref41}, where the probability of a vehicle switching paths is proportional to the relative gap between the best path and the current path. Recently, \cite{ref42} adopted the gap-minimizing method  to the time-dependent UE problem in the transit system with vehicle capacity constraints. This paper formulates a gap-based UE model for the STAP with hard train capacity and interacting OD flows in a multi-line network.

To support this algorithm, the average cost for passengers on a given route and departure from origin are required.
We have the cost $cost(k)$ for passengers whose OD is $k$,
from Eq. (\ref{eq:costr}).

Let us consider the $i^{th}$ passenger whose OD is $k$ on route $r$ and waiting for train $t_1$ on line $l_1 = l^r_{o^k}$ at origin $o^k$.
$q^k(t,r)$ is the number of passengers arriving at their origin to meet train $t$.
The train met by the $i^{th}$ passenger on route $r$ is:
$ t_1 = t^n_{l_1}$ where $n = min(\{j \in T_{l_1} : \sum_{jj \leq j} :q^k(t^{jj}_{l_1},r) \geq i\})$.

Assuming passengers who start earlier also arrive earlier, the train on which the $i^{th}$ passenger whose OD is $k$ on route $r$ arrives at the destination $d^k$ on train $t_2$ and line $l_2 = l^r_{d^k}$ is: \\
$t_2 = t^n_{l_2}$ where $n = min(\{j \in T_{l_2} : \sum_{jj \leq j} :f^k_{pe(d^k, l_2)}(t^n_{l_2}, r) \geq i\})$.

The in-vehicle duration of a route can be calculated by summing the duration of the set of legs on each line on the route.
The waiting cost is the difference between the departure time of the train at the origin, the $i^{th}$ passenger arrived for, and the arrival time at the destination of the train the $i^{th}$ passenger arrived on, minus the in-vehicle time.
The early/late cost can be calculated from the arrival time at the destination.
Thus the cost $c^k_i(t,r)$ of each passenger in the group $q^k(t,r)$ can be calculated based on the output of the simulation.

Finally this enables the average cost $avc^k(t,r)$ for all the passengers in that group $I^k(t, r)$ to be calculated, thus:

Let $avc^k(t, r)$ denote the average cost of option ($t, r$) within OD pair $k$, then it is given by:
\begin{equation}
avc^k(t, r) =  \begin{cases} 
    \frac{\displaystyle \sum_{i \in I^k(t, r)} c_i^k(t, r)}{\displaystyle q^k(t, r)}, \quad \text{\textit{if}} q^k(t, r) > 0\\
    \\
    avc_0^{k}(t, r), \quad  \text{\textit{if}}q^k(t, r) = 0
\end{cases}
\label{eq:29}
\end{equation}
where $avc_k^{0}(t, r)$ is the free flow cost of the option ($t, r$) at $k$, which is calculated when none of the trains on this route are already filled to capacity by passengers on other ODs.
The best option choice within OD pair $k$ is the decision ($t, r$) on line $l_1 = l^r_{o^k}$ with the minimum average travel cost:
\begin{equation}
avc_{*}^k = \min_{t \in T_{l_1}, r \in R^k} avc^k(t, r) 
\label{eq:30}
\end{equation}

The use equilibrium problem on the STAP can be formulated as a mathematical problem (NMP):
\begin{equation}
\min_\mathbf{\textit{Q}}  \sum_{k \in K}\sum_{t \in T_{l_1}} \sum_{r \in R^k}  \left(avc^k(t, r) - avc_{*}^k\right) \times q^k(t, r)
\label{eq:31}
\end{equation}

The objective function in Eq. (\ref{eq:31}) is the total system gap, which is expected to be minimized to find the UE solution. 

\subsection{System Optimum} \label{sec:SO}
Compared to the UE, which focuses on the travel options (time and route)
cost gap, the SO focuses on the total system cost. The objective function of the SO is :
\begin{equation}
\min_\mathbf{\textit{Q}}  \sum_{k \in K} cost(k)
\label{eq:32}
\end{equation}

In general, passengers place high penalties on waiting time, particularly caused by denied boarding under limited service capacities \cite{ref59}. In a system optimal solution, therefore, it is required that all passengers waiting for a train will
be able to board it. The paper aims to find a system to ensure minimal waiting time (with no denied boarding) under the SO solution. Otherwise, passengers can always be shifted to earlier/later departure times given enough supply.


Then, the model presented in Eqs. (\ref{eq:1})-(\ref{eq:costr}) can be simplified if all passengers catch the train they are waiting for.
In this case, Eq. (\ref{eq:chboard}) has its 
inequality replaced by an equation:
\begin{equation}
f^k_{ne(s, l^t)}(t,r) = g^k_s(t, r)
\label{eq:allchboard}
\end{equation}

Constraints in Eqs. (\ref{eq:1})-(\ref{eq:costr}) are updated as follows:

Eqs. (\ref{eq:1}), (\ref{eq:4})-(\ref{eq:7}), and (\ref{eq:9})-(\ref{eq:costr}) keep the same.

Eq. (\ref{eq:2}) is updated to:
\begin{equation}
   g^k_s(t, r) = q^k(t,r)
   \label{eq:24}
\end{equation}

Eq. (\ref{eq:3}) is updated to:
\begin{equation}
g^k_s(t, r) = f^k_{pe(s,l_1)}(t_1,r)
\label{eq:25}
\end{equation}



\section{Evaluation and Solution Procedures} \label{sec:solution-procedure}
\subsection{User Equilibrium}
In this paper, the UE model is formulated with its objective function in Eq. (\ref{eq:31}), subject to the constraints in Eq. (\ref{eq:1})-(\ref{eq:30}). The simulation is used to evaluate the system (simulation is required to determine which passengers catch the train in Eq. (\ref{eq:chboard})). The UE model is solved using an adaptive Frank-Wolfe (AdaFW) algorithm (see Appendix \ref{appendix:customA}). The newly proposed AdaFW algorithm involves two primary loops: \textit{the system-based flow shifting loop} and the \textit{OD-based flow shifting loop}.

\subsection{Approximate System Optimum}
In most existing transit papers, the SO models are solved using simulation-based methods or other heuristic methods.  Therefore, 
these reported `system-optimal' solutions are actually approximations of the exact SO solutions. 
In this paper, we formulate such Approx. SO model with its objective function in Eq. (\ref{eq:32}), subject to the constraints in Eqs. (\ref{eq:1})-(\ref{eq:costr}). We use the simulation to evaluate the system cost, and an AdaFW algorithm to solve the Approx. SO model. A similar approach for obtaining the Approx. SO solution is also described in \cite{ref8}, which proposes a simulation-based SO model solved using a heuristic algorithm.

\subsection{Exact System Optimum}
As mentioned in Section \ref{sec:SO}, this paper assumes the implementation of an SO system without denied boarding. After making this change, \emph{all} the passengers can successfully board the coming train $(dbk_s^k(t) = 0, \quad \forall k \in K, t \in T, s \in S)$, so there is no longer a need to perform a simulation to determine which passenger board and which are denied board and left on the platform. The model is formulated with objective function in Eq. (\ref{eq:32}), subject to the constraints in Eqs. (\ref{eq:1}), (\ref{eq:4})-(\ref{eq:7}), and (\ref{eq:9})-(\ref{eq:costr}), and (\ref{eq:allchboard})-(\ref{eq:25}). The resulting model is now an integer-linear model that can be easily and quickly solved by a standard integer-linear optimization solver. Therefore, the solution is an exact SO solution. 

Broadly speaking, MiniZinc allows problems to be formulated in a way that close to their mathematical formulation \cite{ref43}.
The MiniZinc language lets users write models in a way that is close to a mathematical formulation of the problem, using familiar notation such as existential and universal quantifiers, sums over index sets, or logical connectives like implications \cite{ref57}.
In this paper, the above exact SO model can be solved to optimality by an integer-linear solver. Therefore, the exact SO model is encoded in the MiniZinc modeling language using MiniZinc version 2.8.5 and solved using Gurobi version 9.5.1. \cite{ref47}.

\section{Case study}\label{sec:case-study}
\subsection{Experimental settings}\label{sec:physical-settings}
\begin{figure}[H]
\centering
\includegraphics[width=3.6 in]{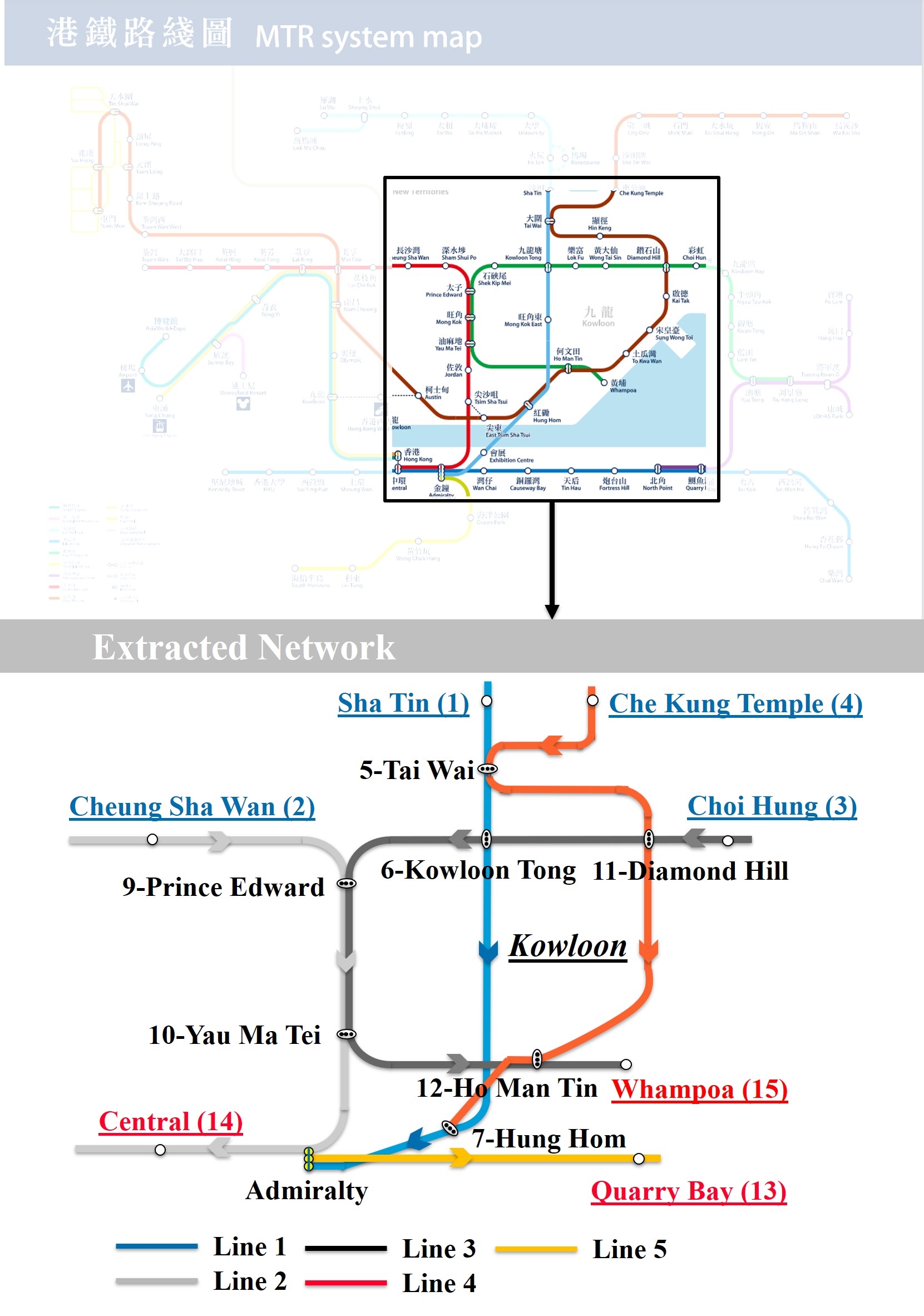}
\caption{Extracted Hong Kong MTR network (adapted from web image: \\ https://www.mtr.com.hk/en/customer/services/system\_map.html)}
\label{Fig4}
\end{figure}

In this paper, we extract a representative subset of the Hong Kong MTR network in the central business district (CBD) with critical links, as shown in Fig. \ref{Fig4}.
 Without loss of generalization, we consider specific OD demands and lines running from north to south, towards the CBD. Other demand flows within the CBD can be incorporated as spatial-temporal background flow to update the corresponding train capacities in our model. Then, the STAP problem can also be solved using the methods presented in this paper.

The modeled section includes 5 lines indicating 26 trains (line 1), 26 trains (line 2), 16 trains (line 3), 35 trains (line 4), and 17 trains (line 5) respectively. Each train has a capacity $S_t =2600$ passengers. The OD matrix (see Table \ref{tab:top_table}) is derived from the smart card data, comprising a total of 52,717 trips. For some OD pairs, a second route will be provided (`D' in Table \ref{tab:top_table}) if two routes have similar transfer requirements for that OD pair, otherwise only a single route is available (`S' in Table \ref{tab:top_table}). Train schedules for \textit{Line 1}-\textit{Line 5} are obtained from the Google Map.

Passengers are expected to make departure time and route decisions
before their trips, including:
\begin{itemize}
\item \textit{Route choice}: Passengers can choose the alternative route for given OD pairs with `\textit{D}' in Table \ref{tab:top_table}.
\item \textit{Departure time choice}: Passengers can choose any train at the origin.
\end{itemize}

The parameter estimates for the individual generalized travel cost functions (\ref{eq:121})-(\ref{eq:123}) are presented in Table \ref{tab:bottom_table}. Parameters $\alpha, \gamma, \mu$ were reported in \cite{ref45}. The parameter $\beta$ is estimated based on attribute weightings where waiting time is weighted as three times the in-vehicle time as described in \cite{ref46}.

\begin{table}[ht]
    \centering
    \begin{minipage}[t]{0.46\textwidth}
        \centering
        \caption{Origin-Destination Matrix(5:50-10:00 AM)}
        \label{tab:top_table}
        \begin{tabularx}{\columnwidth}{cX X X}
            \toprule
            \textbf{O/D} & \textbf{Quarry Bay} & \textbf{Central} & \textbf{Whampoa} \\
            \midrule
            Sha Tin & 5356 (S) & 5663 (D) & 1892 (D) \\
            Cheung Sha Wan & 6116 (S) & 4073 (S) & 2049 (S) \\
            Choi Hung & 14967 (D) & 5852 (S) & 2525 (S) \\
            Che Kung Temple & 1727 (D) & 1848 (D) & 649 (D) \\
            \bottomrule
        \end{tabularx}
    \end{minipage}

    \vspace{0.2cm} 

    \begin{minipage}[b]{0.46\textwidth}
        \centering
        \caption{Parameter values}
        \label{tab:bottom_table}
        \begin{tabularx}{\columnwidth}{cX X}
            \toprule
            \textbf{Attribute} & \textbf{Descriptions}\\
            \midrule
            $\alpha $ & waiting cost per unit time: \$ 18/h\\
            $\beta$ & in-vehicle cost per unit time: \$ 6/h\\
            $\gamma$& early delay cost per unit time: \$ 5/h \\
            $\mu$ & late delay cost per unit time: \$ 12/h\\
    
            \bottomrule
        \end{tabularx}
    \end{minipage}
\end{table}
\subsection{Model comparisons}
\begin{figure}[H]
\centering
\includegraphics[width=3.4 in]{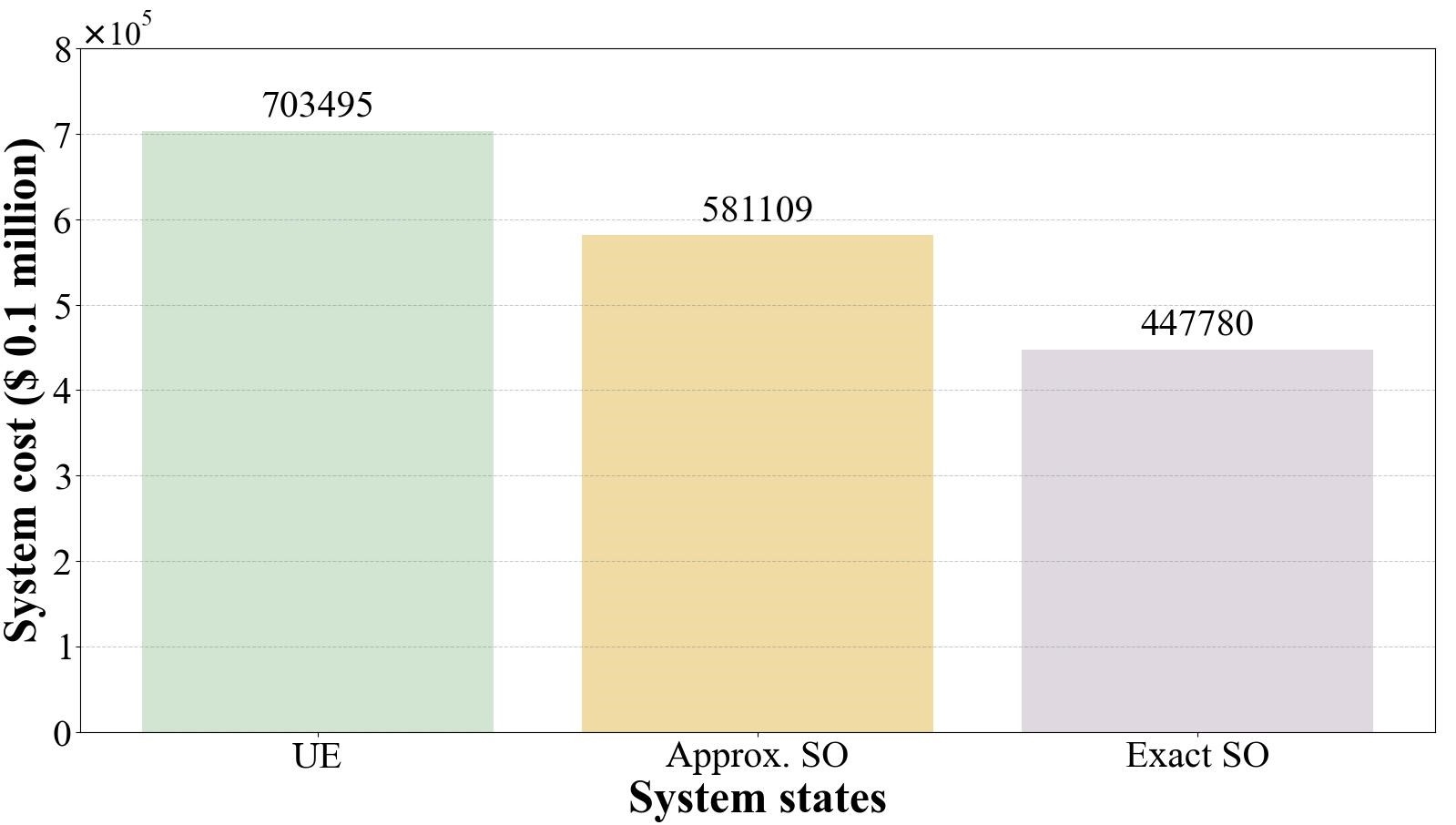}
\caption{Comparison of total system costs among user equilibrium, Approx. SO, and exact SO}
\label{Fig5}
\end{figure}

In this section, we compute a UE solution to determine the congestion cost the system. Then, we compare the UE cost to the costs of Approx. SO and exact SO models. Figure \ref{Fig5} shows that using the exact SO method leads to a minimum system cost of \$447,780, which is 77\% of the cost of the approximate SO solution. When comparing the SO solutions with the UE solution, the exact SO solution shows that the system in the UE state has the potential to reduce its cost by 36.35\%, while the potential improvement is only 17.39\% if the Approx. SO is adopted. Therefore, there is a clear gap in system costs between the exact SO and Approx. SO.
Assume a scenario where a policy maker wants to determine whether providing more flexibility to passengers can improve the system. The flexibility we consider involves the types of available choices for passengers, such as expanding from just departure time choices to include simultaneous departure time and route choices. Currently, we test the departure time choice only case by setting the second routes for all OD pairs as unavailable (updating all OD pairs with `D' in Table \ref{tab:top_table} to `S').

\begin{figure}[H]
\centering
\includegraphics[width=3.4 in]{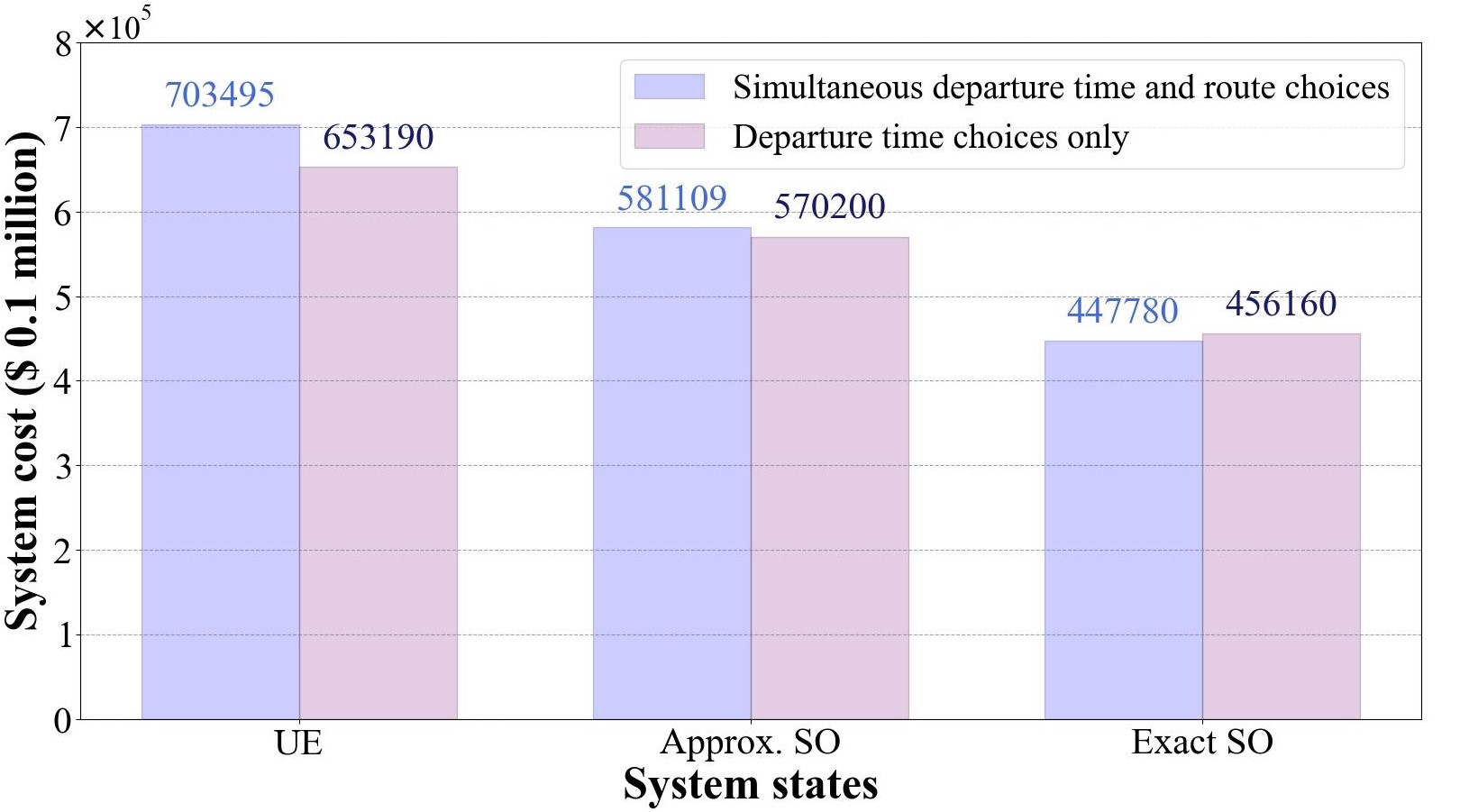}
\caption{System performance comparison: departure time choice only versus simultaneous departure time and route choice}
\label{Fig6}
\end{figure}

From Fig. \ref{Fig6}, we see that the results of both UE and Approx. SO show that more flexibility causes a worse system state (higher system cost indicates worse congestion). However, the exact SO can achieve a better system state with increased flexibility. It is easy to understand the results in UE, as passengers are selfish and make decisions based on their individual benefits. When more flexibility is provided, passengers are expected to make fewer `good' decisions as anticipated by the system. Ideally, in the SO solution, passengers will completely follow the system's suggestions. Thus, the SO solution for the system with only departure time choices is one of the feasible solutions for the system with simultaneous departure time and route choices.
Therefore, if more flexibility is provided, the SO solution is expected to be better (or at least no worse than the SO solution with less flexibility). The exact SO solution aligns with our expectations, while the Approx. SO solution shows the opposite result. This may be caused by the fact that the Approx. SO represents an intermediate state between UE and the exact SO. When the Approx. SO is closer to the UE side, it exhibits similar characteristics to UE. The tendency is also evident in Fig. \ref{Fig6}, where increased flexibility results in a 7.7\% increment in UE compared to only a 1.9\% increment in Approx. SO.

When developing policies for congestion relief, reasonable flexibility should be taken into account. Increased flexibility may enhance the potential of the congested system. However, too much flexibility may complicate the decision-making process, potentially worsening the UE system. 
\subsection{Analysis on targeting passenger groups for shifting}

In this section, we compare the results of UE and SO to identify the most potential passenger group, who will contribute the most to system improvement if they change their decisions as expected by the system.  It is valuable for policy makers to identify which passengers contribute most to congestion. By targeting these passengers, they can propose effective policies to encourage changes in their travel behavior. The main idea of the targeting process is to consider UE as the current state and the SO as the goal state. 
\begin{figure}[H]
\centering
\includegraphics[width=3.5 in]{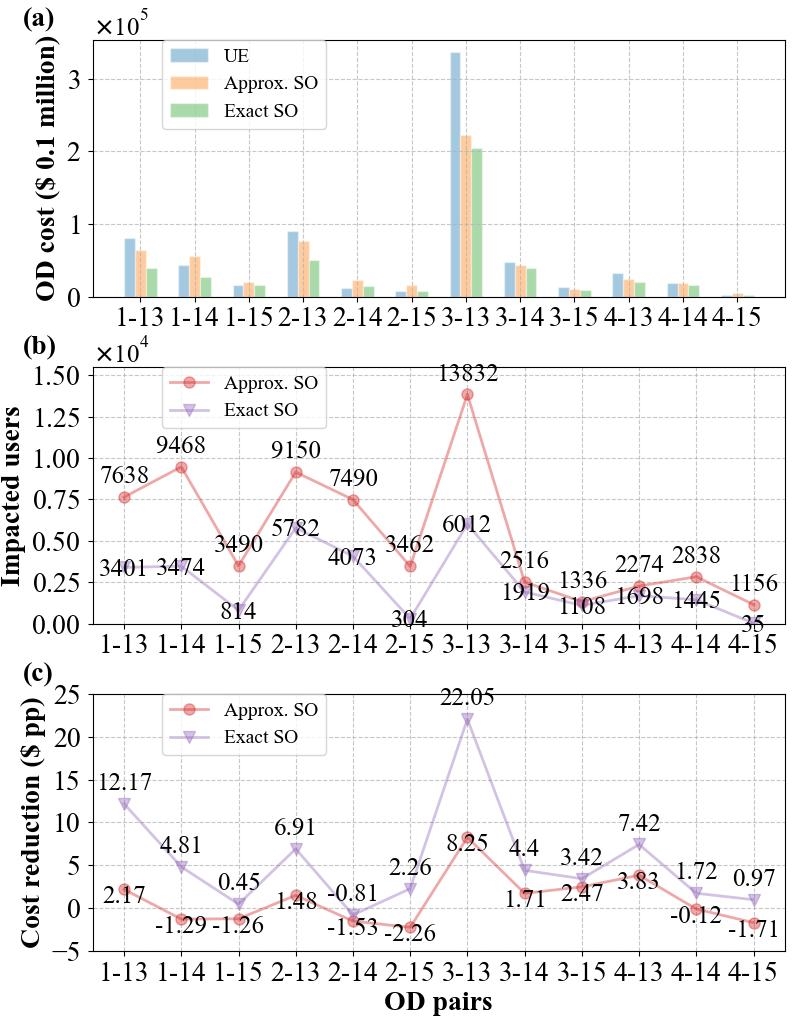}
\caption{Comparison between Approx. SO and exact SO: (a) Travel costs for OD pairs; (b) Number of impacted passengers (potential target passengers for shifting); (c) Cost improvement contribution per passenger}
\label{Fig7}
\end{figure}
Fig. \ref{Fig7}(a) shows the total OD costs for all OD pairs. For simplicity, we assume the corresponding numbers in Fig. \ref{Fig4} represent these OD pairs in the following sections (for example, `3-13' represents the OD pair `Choi Hung-Quarry Bay'). The OD pair `3-13' is found to contribute the most to the system cost, followed by the next top three OD pairs: `2-13', `1-13', and `3-14'. The UE system is expected to reduce the OD cost in most OD pairs, except for OD pair `2-14', where an increase in OD cost is needed compared to the exact SO. However, the Approx. SO provides misleading information, suggesting that the UE system should increase the OD cost in OD pairs  `1-14', `1-15', `2-14', `2-15', and `4-15'. This is opposite to the exact SO, except for the recommendation for `2-14'. Therefore, adopting an inaccurate SO solution may limit the potential improvement of the UE system through intervention. 

Fig. \ref{Fig7}(b) shows the number of impacted passengers by comparing the differences in flow distribution across various options (simultaneous departure time and route choices) between UE and SO. It shows that the top 5 impacted passenger groups using exact SO are located in OD pairs: `3-13' (6012), `2-13' (5782), `2-14' (4073), `1-14' (3474), and `1-13' (3401), while the top 5 impacted passenger groups using Approx. SO are located in OD pairs:`3-13' (13832), `1-14' (9468), `2-13' (9150), `1-13' (7638), and `2-14' (7490). 

To measure the potential cost reduction from each passenger's shift across different OD pairs, we further plot Fig. \ref{Fig7}(c) to show the potential cost contribution per passenger. This is calculated by dividing the potential OD cost gap (obtained from comparing the OD cost difference between UE and SO in Fig. \ref{Fig7}(a) by the corresponding impacted passengers. 

A negative value represents that shifting passengers within that OD pair would increase the corresponding OD cost. The exact SO result shows that passengers in OD `3-13' are the most valuable group for shifting, as a single shift in this group contributes the most to cost reduction (\$22.05) of the system, followed by groups in `1-13' (\$12.17), `4-13' (\$7.42), `2-13' (\$6.91), `1-14' (\$4.81), and `3-14' (\$4.4). The Approx. SO result shows that passengers in OD pair `3-13' are still the most valuable group for shifting with cost reduction (\$8.25) of the system, followed by groups in `4-13' (\$3.83), `3-15' (\$2.47),  `1-13' (\$2.17), `3-14' (\$ 1.71) and `2-13' (\$1.48). There are more negative values in the Approx. SO compared to the exact SO. It is evident that the information from the Approx. SO differs significantly from that of the exact SO, especially when targeting potential passengers. The Approx. SO may suggest targeting less valuable passengers for shifting.

\begin{figure}[H]
\centering
\includegraphics[width=3.4 in]{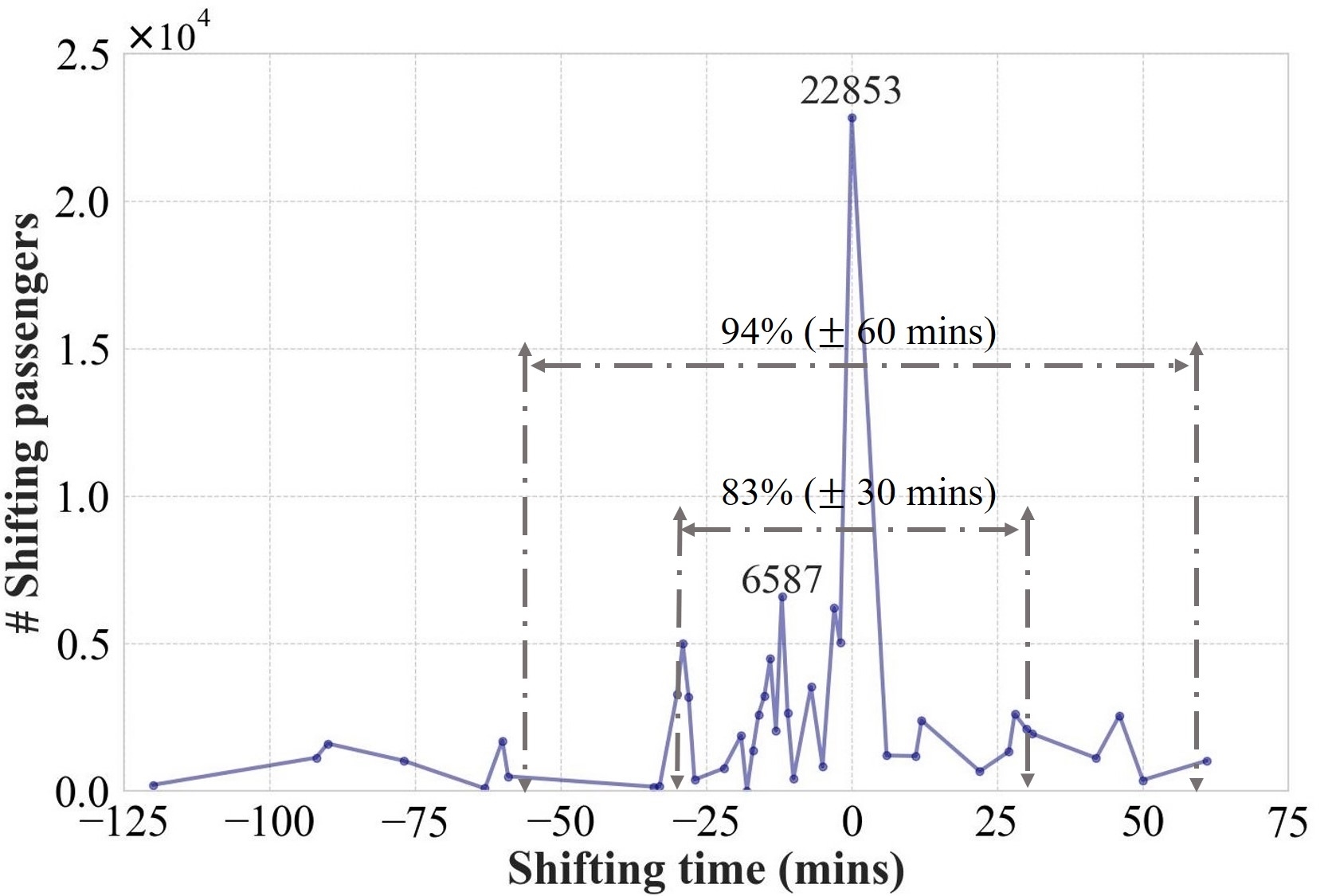}
\caption{The number of passengers needing to be shifted in the system by comparing UE flows and Exact SO flows (negative shifting time (-): shifting to a earlier time; positive (+): shifting to a later time; shifting time (0): no need for shifting)}
\label{Fig12}
\end{figure}

Figure \ref{Fig12} shows the number of passengers needing to be shifted by comparing UE flows with exact SO flows. The result shows that around 80\% passengers (excluding those with a shifting time of 0) need to be shifted within a 30-minute window (-30 to 30 minutes). If we extend the allowable shifting time interval from 30 minutes to 60 minutes (-60 to 60 minutes), then around 94\% passengers (excluding those with a shifting time of 0) can be shifted.

With the above information, policymakers can design individualized incentives rather than generic ones, making incentives more efficient, attractive, and cost-effective.

The Exact SO model can also offer valuable insights into transit operations that involve stop-skipping. OD passenger demand varies significantly along the transit line and throughout the day. For example, busy central stations tend to have a relatively large number of passengers boarding and alighting, while other stations may see fewer passengers \cite{ref60}. Transit scheduling with stop-skipping can reduce passenger travel time and the operational costs for transit service providers \cite{ref61}.

\begin{figure}[ht]
  \centering
  \begin{minipage}[b]{0.48\textwidth}
    \centering
    \includegraphics[width=\textwidth]{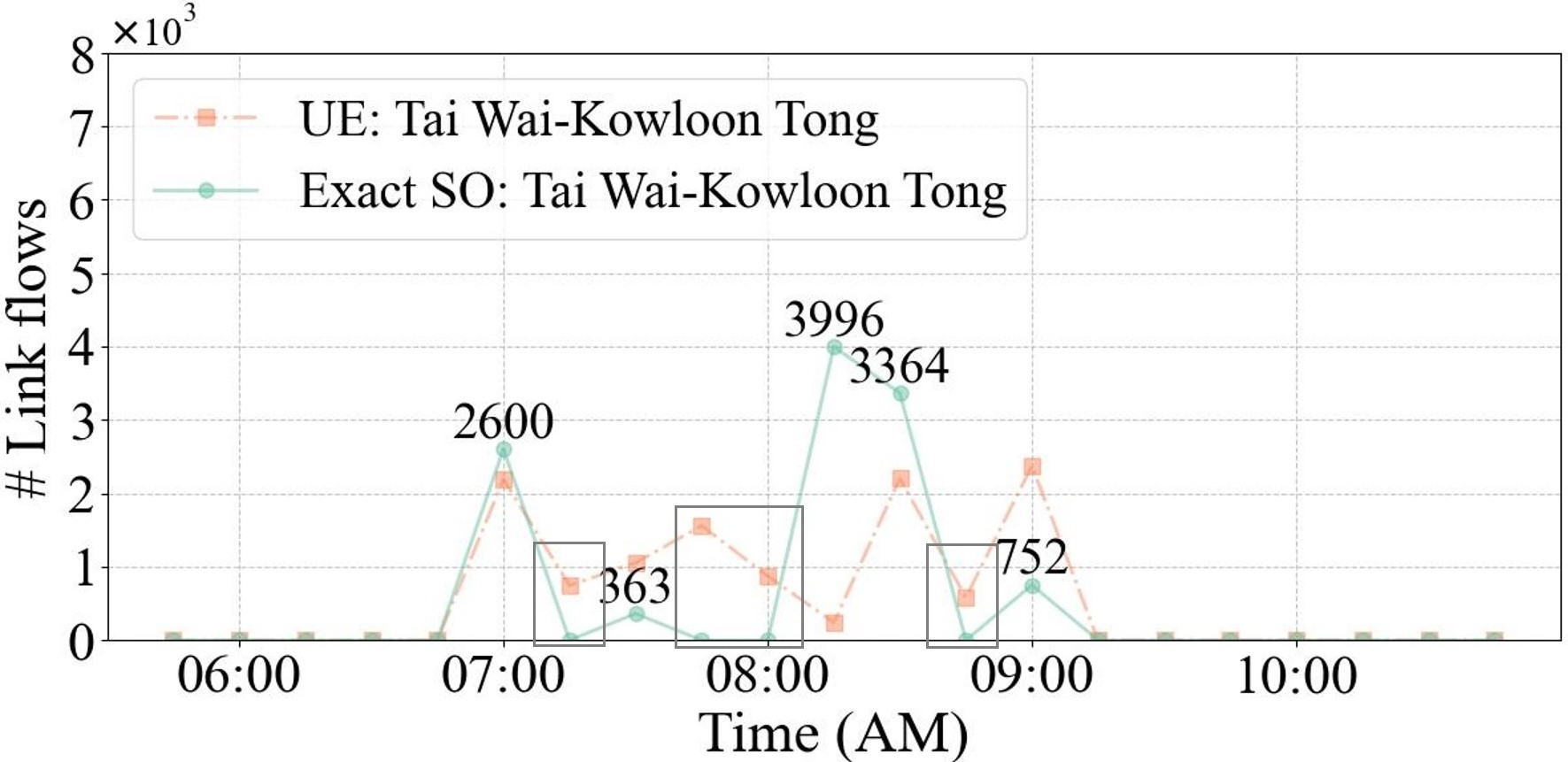} 
    \par\noindent
    \textbf{(a)}
    \label{fig14:minipage1}
  \end{minipage}
  \begin{minipage}[b]{0.48\textwidth}
    \centering
    \includegraphics[width=\textwidth]{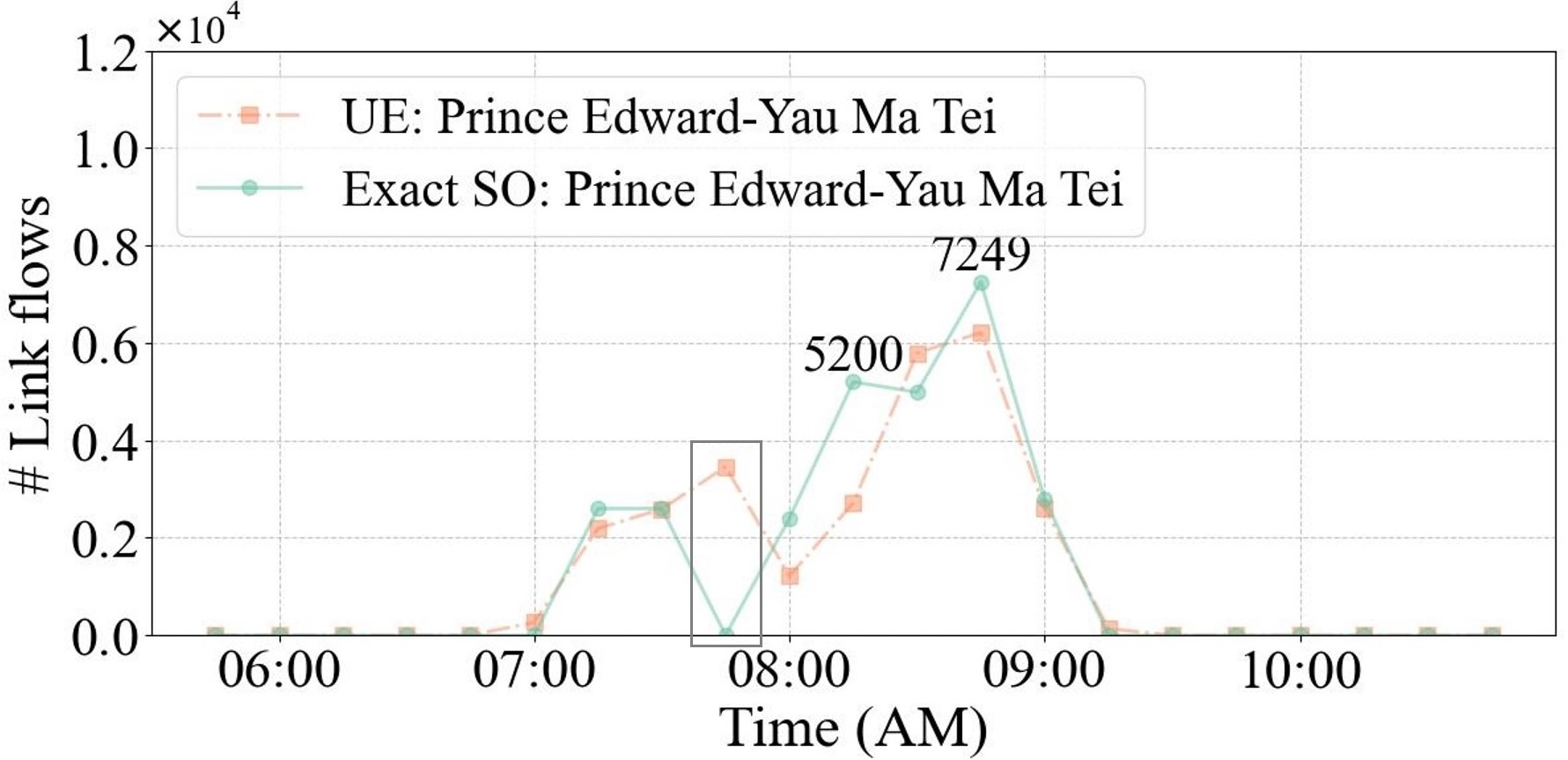} 
    \par\noindent
    \textbf{(b)}
    \label{fig14:main}
  \end{minipage}
  \caption{Link flows over time periods: (a) Link: Tai Wai to Kowloon Tong; (b) Link: Prince Edward to Yau Ma Tei}
  \label{fig14}
\end{figure}

Fig. \ref{fig14} shows link flows over time on critical links: Tai Wai to Kowloon Tong (Fig. \ref{fig14}(a)) and Prince Edward to Yau Ma Tei (Fig. \ref{fig14}(b)). The exact SO solution exhibits significantly different link flows compared to the UE solution. The results show that the exact SO solution indicates certain zero-flow time periods: Tai Wai to Kowloon Tong around 07:15 AM, 07:45 AM, 08:00 AM, and 08:45 AM, and Prince Edward to Yau Ma Tei around 07:45 AM. Trains passing these relevant links during the aforementioned time periods can skip these stops, which can simultaneously improve the system (bringing it closer to the SO solution), reduce the need for incentives for congestion relief, and lower travel costs for downstream passengers.

\subsection{Impact of service capacity}
Apart from the OD demands in Table \ref{tab:top_table}, we also incorporate other demand flows as background flows as mentioned in Subsection \ref{sec:physical-settings}. We simplify these unconsidered demand flows by treating them as fixed flows that occupy a fixed proportion of train capacity. In this section, we measure the potential system performance of UE and exact SO systems at different capacity levels to identify potential improvement areas. The capacity level refers to the actual available capacity after accounting for background flows. A 100\% capacity level indicates that there are no background flows, while a 60\% capacity level means that 40\% of the train capacity (design capacity) is already occupied by background flows. We also consider scenarios where the service provider may upgrade the current trains to higher capacities in the future, examining capacity levels of 120\% and 140\%.
\begin{figure}[ht]
  \centering
  \begin{minipage}[b]{0.48\textwidth}
    \centering
    \includegraphics[width=\textwidth]{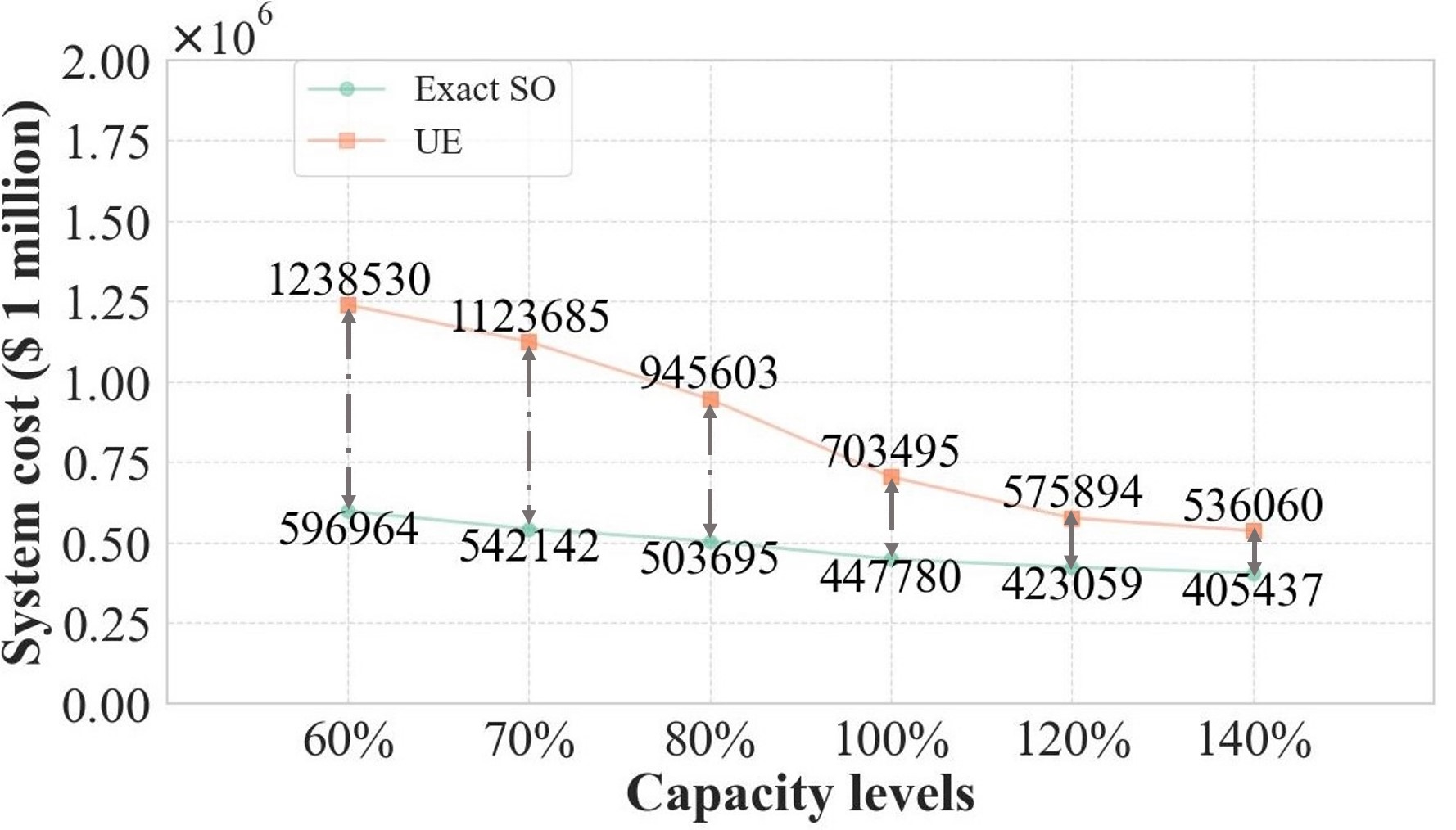} 
    \par\noindent
    \label{Fig30:minipage1}
  \end{minipage}
  \caption{Scenarios on capacity levels 60\%-140\%: System cost evolution curves of UE and Exact SO.}
  \label{Fig30}
\end{figure}

Fig. \ref{Fig30} illustrates the corresponding system costs at different capacity levels. The results indicate that the capacity level has less influence on the SO curve than on the UE curve. This may be due to the fact that, in the SO system, passengers do not experience denied boarding. As a result, higher or lower capacity only affects the early or delay penalties, which are less significant compared to the waiting penalty.  In the UE system, a lower capacity level may result in a significantly higher waiting penalty for passengers. We further evaluate potential system cost improvement ratios (PSCIR) of these UE systems using $1-SO/UE$. The results indicate that UE systems at low capacity levels (i.e., 60\%, 70\%, 80\%) can be significantly improved, with a potential improvement ratio of around 0.5. When the capacity level increases (for example, to 120\% or 140\%), the UE systems can only be improved by around 25\%. This is expected, as the capacity level has less influence on the SO system but a significant impact on the UE system.




\subsection{Impact of demand levels}
In this section, we measure the potential system performance of UE and exact SO systems at varying demand levels to identify potential improvement areas. The demand level refers to the relative demand compared with the current demand: a 100\% demand level refers to the current default demand in Table \ref{tab:top_table}, while a 150\% demand level refers to multiplying the default demand by 1.5. We consider the potential new demand patterns caused by the development of the city, which may stimulate more travel demand in the future. We test these scenarios at 135\%, 150\%, and 165\% demand levels.
\begin{figure}[ht]
  \centering
  \begin{minipage}[b]{0.48\textwidth}
    \centering
    \includegraphics[width=\textwidth]{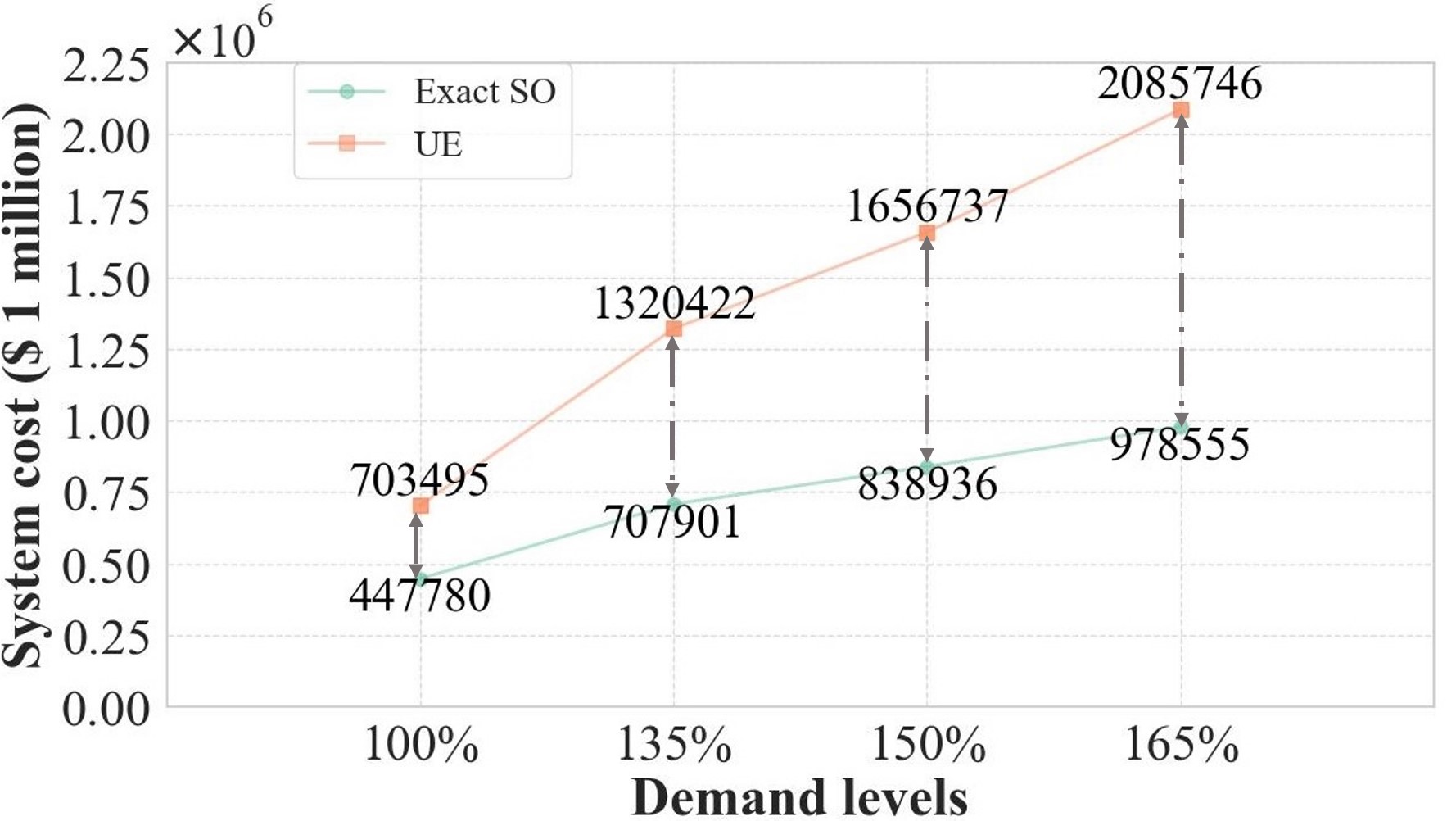} 
    \par\noindent
    \label{Fig8:minipage1}
  \end{minipage}
  \caption{Scenarios on demand levels 100\%-165\%: System cost evolution curves of UE and Exact SO.}
  \label{Fig8}
\end{figure}
Fig. \ref{Fig8} shows a steep, near-linear rise in the UE curve, while the SO curve exhibits a flat, near-linear rise as the demand level increases. This is due to the fact that increased demand leads to more severe competition, complicating the individual decision-making process in UE systems, while it is less influential on the system decision-making process in SO systems. We further evaluate PSCIR of these UE systems: 0.36 (at demand level 100\%), 0.46 (at demand level 135\%), 0.49 (at demand level 150\%), and 0.53 ((at demand level 165\%)). Generally, a higher demand level shows a greater potential for improvement in the UE system to alleviate congestion, with the help of its corresponding SO solution.



We also test scenarios with lower demand levels of 40\%, 60\% and 80\%. The results of these scenarios are discussed in Appendix \ref{appendix:customB}.

\subsection{Model Scalability}
In this section, we implement the exact SO model in scenarios that include all possible additional routes, and increased OD demands. The results are summarized in Tables \ref{tab:bottom_table_2} (include all possible additional routes) and \ref{tab:bottom_table_3} (include all possible additional routes and increased OD demands). 

Considering additional available routes further improves the exact SO solution from \$447,780 to \$406,158, with the running time increasing from 41.395 seconds to 48.613 seconds. Some new routes become popular with high loading in the SO solution, such as `3-6-8-13 (12367)', `3-6-8-14 (2409)', `1-5-11-15 (1727)', and `4-5-8-14 (1569)'. Some new routes contribute little to congestion relief, showing zero loading, such as `1-6-9-8-13 (0)', `4-5-6-15 (0)', and `2-10-15 (0)'. Therefore, adding less valuable routes will not help improve the system. We further test the exact SO model using scenarios that include additional routes at varying demand levels. The results show that the running time increases from 48.613 seconds to 53.038 seconds at a demand level of 135\%, and to 1 minute 4 seconds at a demand level of 150\%. From the perspective of routes and demands, the exact SO demonstrates good scalability in handling more complex scenarios.
\begin{table}[H]
    \centering
    \begin{minipage}[t]{0.46\textwidth}
        \centering
        \caption{Test the exact SO model using the scenario with additional routes}
        \label{tab:bottom_table_2}
        \begin{tabularx}{\columnwidth}{p{0.7cm} p{2.6cm} p{3.8cm}}
            \toprule
            \textbf{Origins} &\textbf{Current Routes}& \textbf{Additional Routes (\& loading)}\\
            \midrule
            \centering \multirow{3}{*}{\textbf{1}} & 1-8-13 & 1-6-9-8-13 (0) \\
            & 1-8-14, 1-6-9-14 & 1-5-11-9-14 (0)\\
            & 1-6-15, 1-5-12-15& 1-5-11-15 (686)\\
            \centering \multirow{3}{*}{\textbf{2}} & 2-8-13 & 2-9-12-7-8-13 (1206)\\
            & 2-14 & \\
            & 2-9-15 & 2-10-15 (0)\\
           \centering \multirow{3}{*}{\textbf{3}} & 3-11-7-8-13, 3-9-8-13 & 3-6-8-13 (12367)\\
            &3-9-14 & 3-6-8-14 (2409), 3-11-7-8-14 (0)\\
            & 3-15&  \\
            \centering \multirow{3}{*}{\textbf{4}} & 4-7-8-13, 4-11-9-8-13 & 4-5-8-13 (1727)\\
            & 4-11-9-14, 4-7-8-14 & 4-5-8-14 (1569)\\
            & 4-12-15, 4-11-15& 4-5-6-15 (0)\\            
            \multicolumn{3}{l}{}\\
            \multicolumn{3}{l}{\textbf{Current exact SO}: System cost \$447,780 (Time: 41.395 seconds)}\\
            \multicolumn{3}{l}{\textbf{Include additional routes} : \$406,158 (Time: 48.613 seconds)}\\
            \multicolumn{3}{l}{Route 1-8-13: 1 (origin) - 8 (transfer) - 13 (destination)}\\
            \bottomrule
        \end{tabularx}
        
    \end{minipage}
        \begin{minipage}[b]{0.46\textwidth}
        \centering
        \caption{Test the exact SO model using scenarios with additional routes at varying demand levels}
        \label{tab:bottom_table_3}
        \begin{tabularx}{\columnwidth}{X X X}
            \toprule
            \textbf{Demand Levels} & \textbf{Running Time} & \textbf{System Cost} \\
            \midrule
            100\% & 48.613 seconds & \$406,158 \\
            135\% & 53.038 seconds & \$636,214\\
            150\% & 1 minute 4 seconds & \$752,755\\           
            \bottomrule
        \end{tabularx}
    \end{minipage}
\end{table}

\section{Conclusion} \label{sec:conclusion}
The schedule-based transit assignment problem (STAP) determines competing flows across departure times and routes using User Equilibrium (UE) or System Optimum (SO) models. A SO model can provide crucial insights into the potential congestion relief of the system, serving as an important reference for further policy design. Current methods for the system-optimum STAP are mainly approximate SO (Approx. SO) models rather than exact SO models, and most of them ignore realistic constraints such as hard capacity and multi-line networks. Unfortunately, these methods may provide misleading information for guiding congestion relief, such as targeting wrong passengers for shifting or underestimating the potential improvement of the system. 

In this paper, we propose an exact SO model for the STAP, accounting for hard train capacity constraints and spatial-temporal competing demand flows in realistic transit networks. A case study of the Hong Kong MTR network is used to validate the proposed method by comparing it with UE and Approx. SO models. We also test scenarios with varying demand and capacity levels using the exact SO model. The model's scalability is further examined by introducing additional routes and increasing OD demands.

Results show that computing an Approx. SO solution for this system indicates a modest potential for congestion reduction measures, with a cost reduction of 17.39\% from the UE solution.  Our exact SO solution is 36.35\% lower than the UE solution, representing more than double the potential for congestion reduction. Compared to the exact SO model, the Approx. SO model underestimates the congestion relief potential of some passengers and targets wrong passengers (less valuable passengers) for shifting. The potential system improvement ratio (1-Exact/UE) increases from 36\% to 53\% as the demand level rises from 100\% (current) to 165\% (future). From the perspective of routes and demands, the exact SO demonstrates good scalability in handling more complex scenarios. The exact SO solution can be used to identify specific opportunities for congestion reduction: (i) which OD pairs have the most potential to gain; (ii) how many passengers can be reasonably shifted; (iii) future system potential given increasing demand and capacity. This information provides insights to help policymakers propose efficient, cost-effective, and individualized incentives for congestion relief.

This work opens the field for several research directions. First, it is of interest to test new policies for congestion relief targeting valuable passengers based on our proposed methods. Second, it is of interest to extend our model to investigate multi-mode transit systems. Third, incorporating the heterogeneity of passengers would be valuable to better reflect real-world features and achieve a more realistic optimal solution.

\section*{Acknowledgments}
We acknowledge Monash University for providing the centrally awarded scholarship. This research is supported by the Australian Research Council under Discovery Project DP190100013. This research is also partly funded by TRENoP (Swedish Strategic Research Area in Transport) Research Center at KTH, Sweden. This research is partially funded by the Australian Government through the Australian Research Council Industrial Transformation Training Centre in Optimisation Technologies, Integrated Methodologies, and Applications (OPTIMA), Project ID IC200100009.




\begin{IEEEbiography}[{\includegraphics[width=1in,height=1.25in,clip,keepaspectratio]{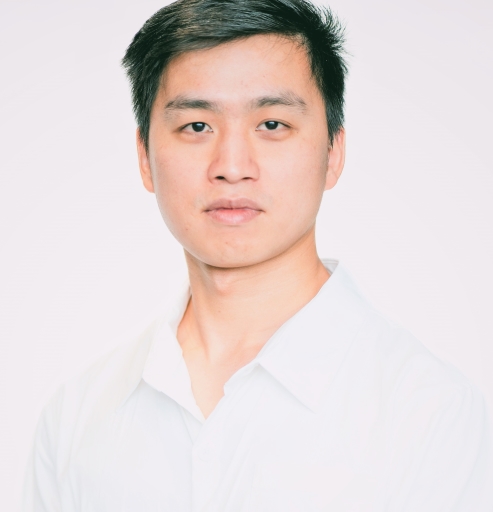}}]{Xia Zhou}
received the B.S. degree in Engineering from Dalian Maritime University, Liaoning, China, and the M.S. degree in engineering from Tianjin University, Tianjin, China. He is pursuing the
Ph.D. degree in the Faculty of Information Technology, Monash University, Melbourne, Australia. His research interests span different techniques and algorithms for optimization, simulation, modeling, and their integration and application to solving transportation assignment problems.\end{IEEEbiography}
\begin{IEEEbiography}[{\includegraphics[width=1in,height=1.25in,clip,keepaspectratio]{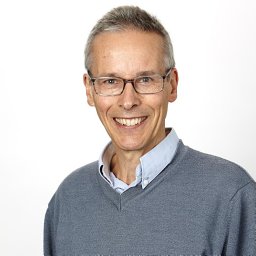}}]{Mark Wallace}
received the bachelor’s degree in mathematics and philosophy from Oxford University, U.K., the master’s degree in artificial intelligence from University of London, U.K., and the
Ph.D. degree from The Southampton University, U.K. He has served as Professor and Associate Dean (Research) with the Faculty ofof Information Technology, Monash University, Melbourne, Australia. His research interests span different techniques and algorithms for optimization and their integration and application to solving complex resource planning and scheduling problems.\end{IEEEbiography}
\begin{IEEEbiography}[{\includegraphics[width=1in,height=1.25in,clip,keepaspectratio]{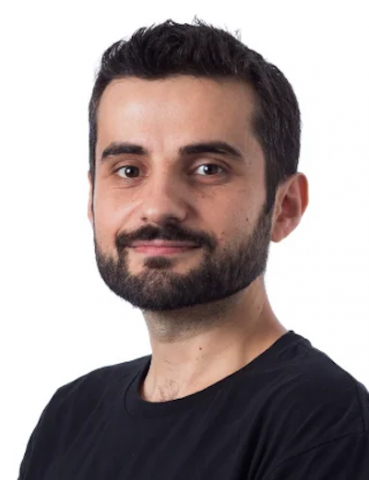}}]{Daniel D. Harabor}
is currently an Associate Professor in the Faculty of Information Technology, Monash University, Melbourne, Australia. He received the B.S. degree (Honours) and the Ph.D. degree in computer science from the Australian National University, Canberra, Australia. His research interests include artificial intelligence, heuristic search, optimisation. In the context of transportation and logistics, he focuses on rail scheduling and capacity analysis of rail supply chains, simulation-based analysis of rail and port operations, data-driven efficacy analysis, and delivery truck fleet optimisation. \end{IEEEbiography}
\begin{IEEEbiography}[{\includegraphics[width=1in,height=1.25in,clip,keepaspectratio]{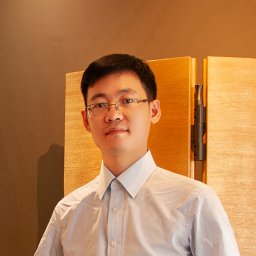}}]{Zhenliang Ma}
received the B.Sc. degree in electrical engineering from Shandong University in 2009, the M.Sc. degree in information technology in 2012, and the Ph.D. degree in transportation engineering from The University of Queensland in 2015. He is currently an Associate Professor of transportation science and a Faculty Member of digital futures with the KTH Royal Institute of Technology. His research interests include statistics, machine learning, computer science-based modeling, simulation, optimization, and control within the framework of selected mobility-related complex systems, which are intelligent transport systems (traffic/public transport/rails) and personal information systems (transport/energy).\end{IEEEbiography}

\begin{appendices}

\section{Formal Constraints}
\label{appendix:customE}
\renewcommand{\theequation}{A.\arabic{equation}} 
\renewcommand{\thefigure}{A.\arabic{figure}}     
\renewcommand{\thetable}{A.\arabic{table}}       
\setcounter{equation}{0} 
\setcounter{figure}{0}   
\setcounter{table}{0}    
This appendix provides formal constraints in Section \ref{sec:constraint}.

The first group of constraints models the enforcement of hard capacity on all legs.

For all lines ($l \in L$), indexes of all trains $i \in T(l)$ on line $l$, and indexes of all stations $j \in S(l)$ on line $l$:
\begin{flalign}
\forall l \in L, &  i \in T(l), j \in S(l): \\
\notag & \sum{k \in K,r \in R^k}: \  f^k_{e^j_l}(t^i_l, r) \leq cap(t) 
\end{flalign}

The second group of constraints models which passengers may wait to board train $t$ at station $s$.

This comprises the passengers in $q^k(t,r)$ whose origin is $o^k$ and the passengers changing trains at $s$. It includes passengers waiting for the first train ($i = 1$), waiting at origin ($s = o^k$), at destination ($s = d^k$), or at other stations ($s \neq o^k \wedge s \neq d^k$).
\begin{flalign}
\forall k \in K, & r \in R^k, l \in L, i \in T(l), j \in S(l): \\ \notag
g^k_{s^j_l}(t^i_l, r)  =& \\ \notag
& (s^j_l=o^k \wedge i=1) \times q^k(t^i_l,r) +\\ \notag
& (s^j_l = o^k \wedge i \neq 1) \times (q^k(t^i_l,r) + db^k_{s^j_l}(t^i_l, r)) + \\ \notag
& (s^j_l=d^k) \times 0 + \\ \notag
& (s^j_l \neq o^k \wedge s^j_l \neq d^k) \times \\ \notag 
&\sum_{l_1 \in L, x \in T(l_1), y \in S(l_1)} : f^k_{e^{y-1}_{l_1}}(t^x_{l_1}, r) \times \\ \notag
& (ch(r,s^j_l,l_1,l) \wedge con(s, t^x_{l_1},t^i_l) \wedge s^j_l=s^y_{l_1})) 
\end{flalign}

The third group of constraints models which passengers may be left behind and must wait for train $t$ at station $s$. 

It includes the denied boarding passengers for the first train ($i = 1$), after the last train ($i = tct_l + 1$) and other trains ($i \neq 1 \wedge i \neq tct_l + 1$).
\begin{flalign}
\forall k \in K, & r \in R^k, l \in L, i \in T(l), j \in S(l): \\ \notag
db^k_{s^j_l}(t^i_l, r)  =& \\ \notag
(i = 1 & \wedge i = tct_l + 1) \times 0 +\\  \notag
(i \neq 1 & \wedge i \neq tct_l + 1) \times \left(g^k_{s_l^j}(t_l^{i-1}, r) - f^k_{e^j_i}(t^{i-1}, r)\right) 
\end{flalign}

The fourth group of constraints models the passengers on leg $e$ aboard train $t$. 

For station $s$ and route $r$, we use $tra(s, r)$ as a Boolean algebra to check the relationship between $s$ and $r$. $tra(s, r)$ is true if $s$ is a transfer station on route $r$. Passengers on the exiting edge are the same as that on entering edge, for station $s$ which is not a  transfer station on route $r$ ($tra(s, r) = False$), or $s$ is not a origin ($j \neq 1$, $s \neq o^k$), or destination of  ($s \neq s^{sct_l}_l$, $s \neq d^k$) line $l$ or OD $k$. For other stations, the passengers ($ simul_{ne(s, l)} ^k(t, r)$) exiting edge $e$ are decided by the simulation. 
\begin{flalign}
 \forall k \in K,& r \in R^k, l \in L, i \in T(l), j \in S(l): \\ \notag
 f_{ne(s^j_l, l)}^k(t^i_l, r) = & \\ \notag 
  f^k_{pe(s^j_l, l)}& (t^i_l, r) \times (j \neq sct_l \wedge j \neq 1 \wedge s^j_l \neq d^k \wedge  \\ \notag 
  &  tra(s^j_l, r) = False \wedge s^j_l \neq o^k) + \\ \notag
  simul_{ne(s^j_l, l)} ^k&(t^i_l, r) \times (j = sct_l \lor j = 1 \lor s^j_l = d^k \lor \\ \notag
  & tra(s^j_l, r) \lor s^j_l = o^k) \wedge \\ \notag
 & (simul_{ne(s^j_l, l)}^k(t^i_l, r) \in \mathbb{Z}^{0+} \wedge \\ \notag
 & simul_{ne(s^j_l, l)}^k(t^i_l, r) \leq g^k_{s^j_l}(t^i_l, r))
\end{flalign}

The fifth group of constraints models the passengers on legs exiting origins. 

All OD demands $\{q^k, k \in K\}$ are assigned to trains ($i \in T(l)$) with routes ($r \in R$).

\begin{flalign}
 \forall k \in K &, l \in \{l^r_{o^k}\}: \\ \notag
 q^k = & \sum r \in R^k, i \in T(l):  q^k(t^i_l, r) = \\ \notag
  &\sum r \in R^k, i \in T(l): 
   f^k_{ne(o^k, l)}(t^i_l, r)
\end{flalign}

The sixth group of constraints models the passengers on legs entering destinations. 

All arrival passengers at destinations from entering leg $e$ are consistent to the corresponding OD demands.
\begin{flalign}
 \forall k \in K & : \\ \notag
  q^k = & \sum r \in R^k, l_1 \in \{l^r_{o^k}\}, i_1 \in T(l_1):  q^k(t^{i_1}_{l_1}, r) = \\ \notag
  & \sum r \in R^k, l_2 \in \{l^r_{d^k}\}, i_2 \in T(l_2): f^k_{pe(d^k, l_2)}(t^{i_2}_{l_2}, r)
\end{flalign}

The seventh group of constraints models summed passengers on leg $e$, and denied boarding passengers at station $s$.  

For passengers on leg $e$:
\begin{flalign*}
  \forall k \in K,& l \in L, i \in T(l), j \in S(l): \\ \notag
  & fk^k_{e^j_l}(t^i_l) = \sum r \in R^k:  f^k_{e^j_l}(t^i_l, r)
\end{flalign*}

For denied boarding passengers at station $s$:
\begin{flalign*}
  \forall k \in K,& l \in L, i \in T(l), j \in S(l): \\ \notag
  & dbk^k_{s^j_l}(t^i_l) = \sum r \in R^k:  db^k_{s^j_l}(t^i_l, r)
\end{flalign*}

The eighth group of constraints models total in-vehicle cost of passengers whose OD is $k$.  

For in-vehicle cost on leg $e$ on train $t$:
\begin{flalign*}
  \forall k \in K,& l \in L, i \in T(l), j \in S(l): \\ \notag
  & ivcst_{e^j_l}(t^i_l) = \alpha \times fk^k_{e^j_l}(t^i_l) \times dur(e^j_l)
\end{flalign*}

For total in-vehicle cost of all passengers whose OD is $k$:
\begin{flalign}
  \forall k \in K &: \\ \notag
  & ivcst(k) = \sum l \in L, i \in T(l), j \in S(l): ivcst_{e^j_l}(t^i_l)
\end{flalign}

The ninth group of constraints models  total waiting cost of passengers whose OD is $k$.  

For waiting cost caused by time gap between two close connecting trains (from train $t_1$ to train $t_2$, where $t_2$ is the earliest train connecting to $t_1$) at station $s$.
\begin{flalign}
  \forall k \in K,  l \in L&, i \in T(l), j \in S(l): \\ \notag
    \forall l_2 \in L, & i_2 \in T(l_2): chcst^k_{s^j_l}(t^i_l, t^{i_2}_{l_2}) = \\ \notag 
 \sum r \in  &  R^k: f^k_{pe(s^j_l, l)}(t^i_l, r) \times \beta \times \\ \notag
  & (dep(t^{i_2}_{l_2}, s^j_l) - arr(t^i_l, s^j_l))  \times \\ \notag
  & (ch(s^j_l, r, l, l_2) \wedge con(s^j_l, t^i_l, t^{i_2}_{l_2}))
\end{flalign}

For total waiting cost caused by time gap between two close connecting trains for all passengers whose OD is $k$:
\begin{flalign}
    \forall k \in K &: \\ \notag
     chcst& (k) = \sum s \in S, t_1, t_2 \in T: \\ \notag
    & chcst1^k_s(t_1, t_2) \times ch(s, r, l^{t_1}, l^{t_2})
\end{flalign}

For total waiting cost caused by denied boarding by train $t$ for all passengers whose OD is $k$:
\begin{flalign}
    \forall k \in K: \\ \notag
     dbcst&(k)  = \\ \notag
    & \sum l \in L, n \in S_l: (tt^n_l - tt^{n-1}_l) \times \\ \notag
     & dbk^k_s(t^n_l) \times \beta
\end{flalign}

The tenth group of constraints models total early/late penalty cost for all passengers when they arrive at their destinations. 

For early/late penalty cost on line $l$ at station $s$:
\begin{flalign}
    \forall k \in K , l \in L , n \in S_l & : \\ \notag
    elcst1^k(l, n) = & fk^k_{pe(d^k, l)}(t^n_l) \times \\ \notag
    & (\gamma \times \max(540 - arr(t^n_l, d^k), 0) + \\ \notag
    & \mu \times \max (arr(t^n_l, d^k) - 540, 0))
\end{flalign}

For total early/late penalty cost for all passengers whose OD is $k$:
\begin{flalign}
    \forall k \in K: & \\ \notag
    & elcst(k) = \sum l \in L, n \in T_l: elcst1^k(l, n) 
\end{flalign}

\begin{flalign}
    \forall k \in K: & \\ \notag
    cost&(k) = ivcst(k) + chcst(k) + dbcst(k) + elcst(k)
\end{flalign}

\section{Solution algorithm for UE and Approx. SO}
\label{appendix:customA}

\renewcommand{\theequation}{B.\arabic{equation}} 
\renewcommand{\thefigure}{B.\arabic{figure}}     
\renewcommand{\thetable}{B.\arabic{table}}       
\setcounter{equation}{0} 
\setcounter{figure}{0}   
\setcounter{table}{0}    

This appendix provides additional details on our proposed the adaptive Frank-Wolfe (AdaFW) solution algorithm for the UE and approximate SO models. 

In this paper, we design an AdaFW algorithm to solve both the UE and Approx. SO problems. The newly proposed AdaFW algorithm involves two primary loops: \textit{the system-based flow shifting loop} and the \textit{OD-based flow shifting loop}. The algorithm will first run the system-based flow shifting loop until the convergence criterion is met, and then use its output as the input for the OD-based flow shifting loop. The output of the OD-based flow shifting loop is the final solution. 
\begin{algorithm}
\caption{AdaFW algorithm for UE and Approx. SO} 
\label{alg:solutionalgorithm}
\begin{algorithmic}[1]
\STATE \textbf{Procedure}{ (1) System-based flow shifting loop:}
\STATE {\textsc{Input}}: $\textit{Q}$
\STATE {\textsc{Output}}: $\mathbf{\textit{Q}_\textit{1}^*}$\\
\FOR{$j \leftarrow 0$ to $N$}
\STATE {\textsc{Evaluation}}: 
\STATE \hspace{0.5cm}$\mathbf{\textit{F}_\textit{j}} \gets \mathbf{\textit{Q}_\textit{j}}$, via simulation; 
\STATE \hspace{0.5cm}$\mathbf{\textit{C}_\textit{j}} \gets \mathbf{\textit{F}_\textit{j}}$, via Eqs. (\ref{eq:11})-(\ref{eq:costr}); 
\STATE \hspace{0.5cm}$\mathbf{\zeta(\mathbf{\textit{Q}_\textit{j}})} \gets \mathbf{\textit{C}_\textit{j}}$, \textit{Objective function:} via Eq. (\ref{eq:31}) (UE) or Eq. (\ref{eq:32}) (Approx. SO)

\STATE {\textsc{Direction finding}}: $\mathbf{\textit{C}_\textit{j}^*} \gets \mathbf{\textit{C}_\textit{j}} $, via Eq. (\ref{eq:30})
\STATE \hspace{0.5cm}{$\mathbf{\textit{V}_\textit{j}} \gets \mathbf{\textit{C}_\textit{j}^*}$}, shifting all passengers from the bad to the best option choices

\STATE {\textsc{Step-size determination}}: 
\STATE \hspace{0.5cm}$\mathbf{\theta_\textit{j}^*} $, via Golden section search method \cite{ref44}; 
\STATE \hspace{0.5cm}$\mathbf{\sigma_\textit{j}} \gets \mathbf{\theta_\textit{j}^*}$ via Eqs. (\ref{eq:49})-(\ref{eq:52})
\STATE {\textsc{Update}}: $\mathbf{\textit{Q}_\textit{j+1}} = \mathbf{\textit{Q}_\textit{j} + \mathbf{\sigma_\textit{j}}(\mathbf{\textit{V}_\textit{j}} - \mathbf{\textit{Q}_\textit{j}})}$;

\STATE \hspace{0.5cm}$\mathbf{\zeta(\mathbf{\textit{Q}_\textit{j+1}})}\gets \mathbf{\textit{Q}_\textit{j+1}}$, repeat \textsc{Evaluation} step (Line 5)
\STATE {\textsc{Convergence test}}:
\IF{$\mathbf{\zeta(\mathbf{\textit{Q}_\textit{j+1}})} \leq \mathbf{\zeta(\mathbf{\textit{Q}_\textit{j}})}$} 
    \STATE $\mathbf{\mathbf{\textit{Q}_\textit{1}^\textit{*}}}=\mathbf{\textit{Q}_\textit{j+1}}$, 
        $j = j+1$, 
        $\mathbf{\textit{Q}_\textit{j}} = \mathbf{\textit{Q}_\textit{1}^*}$
    \ELSE 
        \STATE{break}
    \ENDIF

\ENDFOR
\RETURN $\mathbf{\textit{Q}_\textit{1}^*}$
\STATE \textbf{End Procedure} {(1)}
\STATE{}
\STATE \textbf{Procedure}{ (2) OD-based flow shifting loop: }
\STATE {\textsc{Input}}: $\textit{Q}_\textit{1}^*$ from \textbf{Procedure} (1) 
\STATE {\textsc{Output}}: $ \mathbf{\textit{Q}^*}$, $\mathbf{\zeta(\mathbf{\textit{Q}^\textit{*}})}$\\
\FOR {$i \gets 0$ to $N$}
\STATE {\textsc{Evaluation}}: $\mathbf{\textit{F}_\textit{i}} \gets \mathbf{\textit{Q}_\textit{i}} $; $\mathbf{\textit{C}_\textit{i}} \gets \mathbf{\textit{F}_\textit{i}}$;
\STATE \hspace{0.5cm}  $\mathbf{\zeta(\mathbf{\textit{Q}_\textit{i}})} \gets \mathbf{\textit{C}_\textit{i}}$
\STATE {\textsc{select randomly:}} $k \in O \times D$
\STATE {\textsc{Direction finding}}: (\textit{only within OD pair $k$})
\STATE \hspace{0.5cm}$\mathbf{\textit{C}_\textit{i}^*(\textit{k})} \gets \mathbf{\textit{C}_\textit{i}(\textit{k})} $, {$\mathbf{\textit{V}_\textit{i}(\textit{k})} \gets \mathbf{\textit{C}_\textit{i}^*(\textit{k})}$}
\STATE {\textsc{Step-size determination}}: (\textit{only within $k$})
\STATE \hspace{0.5cm}$\mathbf{\theta_\textit{i}^*(\textit{k})} $ via Golden section search method; 
\STATE \hspace{0.5cm}$\mathbf{\sigma_\textit{i}(\textit{k})} \gets \mathbf{\theta_\textit{i}^*(\textit{k})}$ via Eqs. (\ref{eq:49}), (\ref{eq:52})-(\ref{eq:53})
\STATE{\textsc{Update}}: (\textit{only update flows within OD $k$})
\STATE \hspace{0.5cm}$\mathbf{\textit{Q}_\textit{i+1}} = \mathbf{\textit{Q}_\textit{i} + \mathbf{\sigma_\textit{i}(\textit{k})}(\mathbf{\textit{V}_\textit{i}} - \mathbf{\textit{Q}_\textit{i}})}$
\STATE \hspace{0.5cm}$\mathbf{\zeta(\mathbf{\textit{Q}_\textit{j+1}})} \gets \mathbf{\textit{Q}_\textit{i+1}}$, repeat the \textsc{Evaluation} (Line 30)
\STATE {\textsc{Convergence test}}:
    \IF{$\mathbf{\zeta(\mathbf{\textit{Q}_\textit{i+1}})} \leq \mathbf{\zeta(\mathbf{\textit{Q}_\textit{i}})}$ } 
    \STATE $\mathbf{\mathbf{\textit{Q}_\textit{2}^\textit{*}}}=\mathbf{\textit{Q}_\textit{i+1}}$, 
        $i = i+1$, 
        $\mathbf{\textit{Q}_\textit{i}} = \mathbf{\textit{Q}_\textit{2}^*}$
    \ELSIF{$i \leq N_2$}
    \STATE {$i = i+1$, $\mathbf{\textit{Q}_\textit{i}} = \mathbf{\textit{Q}_\textit{2}^*}$}
    \ELSE 
        \STATE{break}
    \ENDIF
\ENDFOR
\RETURN $ \mathbf{\textit{Q}_\textit{2}^\textit{*}}$, $\mathbf{\zeta(\mathbf{\textit{Q}_\textit{2}^\textit{*}})}$, $\mathbf{\textit{Q}^\textit{*}} = \mathbf{\textit{Q}_\textit{2}^*}$, $\mathbf{\zeta(\mathbf{\textit{Q}^\textit{*}})} = \mathbf{\zeta(\mathbf{\textit{Q}_\textit{2}^\textit{*}})}$
\STATE \textbf{End Procedure} {(2)}
\STATE {($N$, $N_2$ are large constant integers)}
\end{algorithmic}
\end{algorithm}
In this algorithm, the following variables are used:
\ \\

$\mathbf{\textit{Q}}$ is the passenger demand distribution over trains and  routes at origins, $\mathbf{\textit{Q}} = \{q^k(t, r): k\in K, r \in R^k, t \in T_l$, where $l = l^r_{o^k}\}$;
\ \\

$\mathbf{\textit{F}}$ is the passenger flow distribution over legs connecting to destinations, $\mathbf{\textit{F}} = \{f^k_{pe(d^k, l)}(t, r): k \in K, r \in R^k, t \in T_l$, where $l = l^r_{d^k}\}$;
\ \\

$\mathbf{\textit{C}}$ is the travel cost distribution over:
\begin{itemize}
    \item trains and routes at origins for UE, $\mathbf{\textit{C}} =\{avc^k(t, r): k\in K, r \in R^k, t \in T_l$, where $l = l^r_{o^k}\}$;
    \item OD pairs for Approx. SO, $\mathbf{\textit{C}} = \{cost(k): k \in K\}$;
\end{itemize}
\ \\
$\mathbf{\zeta(\mathbf{\textit{Q}})}$ is the value of objective function with the input variable $\mathbf{\textit{Q}}$:
\begin{itemize}
    \item total system cost gap for UE: via Eq. (\ref{eq:31});
    \item total system cost for Approx. SO: via Eq. (\ref{eq:32});
\end{itemize}

From the perspective of system-based flow shifting loop, it follows an `all-at-once' approach. It simultaneously updates solutions for all OD pairs using:
\begin{itemize}
    \item the system step size $\theta$;
    \item OD weighted ratio distribution $\mathbf{\partial} = \{\partial^k: k \in K\}$;
    \item weighted time option ratio distribution \
    
    $\boldsymbol\varphi = \{\varphi^{k}(t, r): k \in K, (t, r) \in \mathbf{\tau}^k\}$
\end{itemize}
to update the actual step size $\mathbf{\sigma} = \{\sigma^{k}(t, r):  k \in K, r \in R^k, t \in T_l$, where $l = l^r_{o^k}\}$. $\mathbf{\tau}^k$ is the set of non-optimal options including simultaneous departure time and route choices. These $\sigma^{k}(t, r)$, $\partial^k$ and $\varphi^{k}(t, r)$ are updated by Eqs. (\ref{eq:49})-(\ref{eq:52}).
\begin{equation}
\label{eq:49}
\sigma^{k}(t, r) = \theta * \partial^k * \varphi^{k}(t, r)
\end{equation}
This calculation involves the step size $\theta$ (from the system perspective to determine the optimal system step size), OD relative gap $\partial_k$ (the gap relative to the `ideal' UE for OD pair $(k, t)$, where all used time options have the same and lowest option cost), and the time option ratio $\varphi^{k}(t, r)$ (to determine how to target passengers in non-optimal options for shifting). 
\begin{equation}
\label{eq:50}
C^k(t_{\tau}, r_{\tau}) = \frac{ \displaystyle \sum_{(t, r) \in \mathbf{\tau}^k} avc^k(t, r)}{n}
\end{equation}
where $C^k(t_{\tau}, r_{\tau})$ is the average option cost of all non-optimal options $\mathbf{\tau}^k$. $n$ represents the total number of available options within OD pair $k$.
\begin{equation}
\label{eq:51}
\partial^k = \frac{\displaystyle C^k(t_{\tau}, r_{\tau})- avc^k_*}{C^k(t_{\tau}, r_{\tau})}
\end{equation}
In Eq. (\ref{eq:51}), $avc^k_*$ is the minimum option cost within OD pair $k$. This is used to measure the gap between the average cost and the optimal cost within OD pair $k$. For any given value of $avc^k_*$, $\partial_k$ is minimised by minimising $C^k(\tau_k)$. A larger gap will stimulate a greater proportion of passengers to be shifted within this OD pair $k$.
\begin{equation}
\label{eq:52}
\varphi_{k}(t, r) = \frac{avc^k(t, r)}{{\displaystyle \sum_{t^{'} \in T_{l}, l = l^r_{o^k}} \sum_{r^{'}\in R_k} avc^k(t^{'}, r^{'})}}
\end{equation}
Eq. (\ref{eq:52}) indicates that the number of passengers from non-optimal options waiting to be shifted depends on the weighted ratio of the time options within OD pair $k$. A non-optimal time option with a higher option cost will result in a greater number of passengers being shifted to the optimal option.

From the perspective of the OD-based flow shifting loop, it represents an alternative flow update strategy known as `one-at-a-time,' where solutions are updated for a single OD pair per iteration. Most steps of this loop are the same as those in the system-based flow shifting loop, except for the determination of the step size. Eq. (\ref{eq:51}) is replaced with:
\begin{equation}
\partial^k =  \begin{cases} 
    1, \quad \text{if } k \text{ is targeted}\\
    \\
    0, \quad  \text{if } k \text{ is not targeted}
\end{cases}
\label{eq:53}
\end{equation}
The OD-based flow shifting loop ends when transferring just one passenger from a non-optimal option (departure time and route) to any optimal option, on any OD, fails to improve the value of the objective function. 

 




\section{Effects of demand levels 40\%, 60\%, and 80\%}
\label{appendix:customB}
\renewcommand{\theequation}{C.\arabic{equation}} 
\renewcommand{\thefigure}{C.\arabic{figure}}     
\renewcommand{\thetable}{C.\arabic{table}}       
\setcounter{equation}{0} 
\setcounter{figure}{0}   
\setcounter{table}{0}    

This appendix provides extra case studies on  demand levels 40\%, 60\%, and 80\%.

\begin{figure}[H]
\centering
\includegraphics[width=3.6 in]{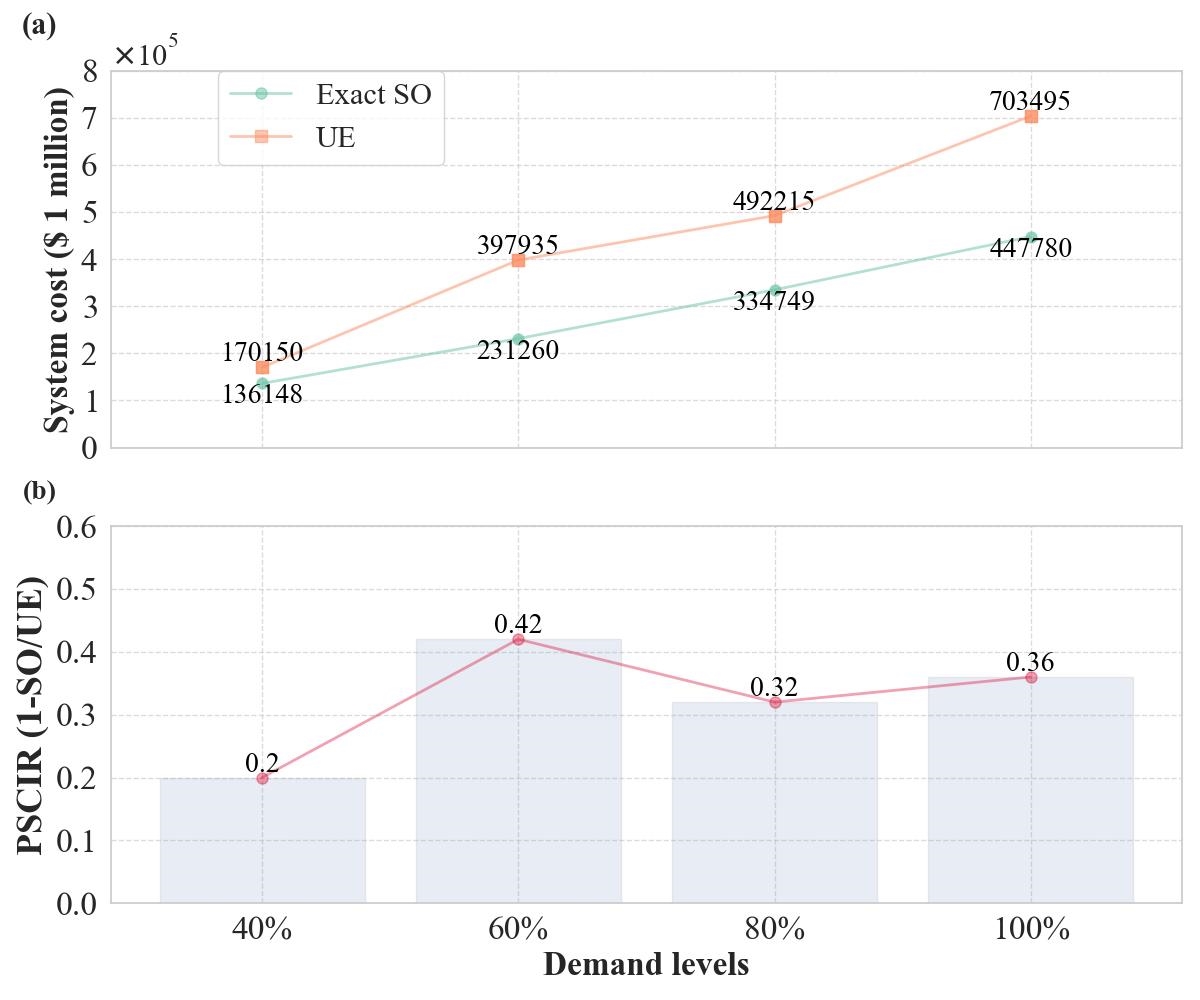}
\caption{Scenarios on demand levels 40\%-100\%: (a) System cost evolution curves of UE and Exact SO; (b) Potential system cost improvement ratios (PSCIR) for UE systems}
\label{Fig20}
\end{figure}

\begin{figure}[H]
\centering
\includegraphics[width=3.4 in]{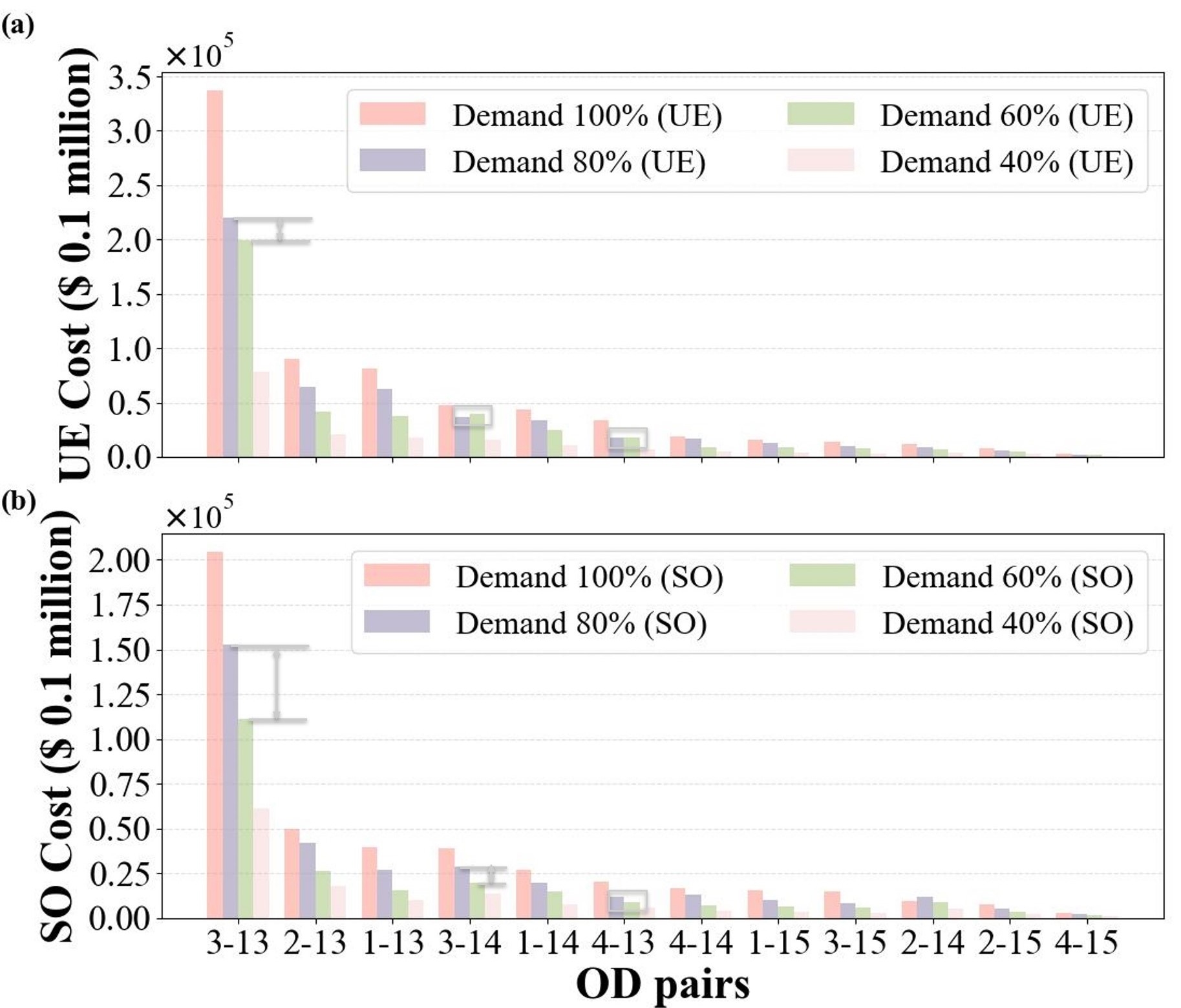}
\caption{OD performance in UE and SO systems across 4 demand levels: (a) Travel cost for each OD pair (UE); (b) Travel cost for each OD pair (SO)}
\label{Fig9}
\end{figure}

As shown in Fig. \ref{Fig20}, there is a unique scenario at the 60\% demand level, characterized by a high potential improvement of 0.42. To explore the scenario at the 60\% demand level, we evaluate OD costs of its surrounding scenarios at 100\%, 80\%, 60\%, and 40\% demand levels, as shown in Fig. \ref{Fig9}. 

Figs. \ref{Fig9}(a) and (b) show that OD costs for almost all OD pairs decrease with the reduction in demand levels in both UE and SO systems, although the reduction ratios vary significantly across OD pairs. In the UE system, when reducing the demand level from 80\% to 60\%, the cost reduction for OD pairs `3-13', `3-14', and `4-13' at the 60\% demand level is not obvious. This indicates that reducing the demand level from 80\% to 60\% has a slight impact on the cost reduction for these OD pairs. 

However, in the SO system, when reducing the demand level from 80\% to 60\%, the cost reduction for OD pairs `3-13', `3-14', and `4-13' is significant. Conversely, at both the 80\% and 40\% demand levels, reducing from a higher demand level results in significant cost reductions in both UE and SO systems. Since the summed OD cost of `3-13', `3-14', and `4-13' accounts for nearly 58\% of the total system cost at the 60\% demand level, a high OD cost in the UE system and a low OD cost in the SO system lead to the high potential improvement (0.42) at the 60\% demand level, as shown in Fig. \ref{Fig8}(b), compared to the 80\% and 40\% demand levels (0.32 and 0.20).

\begin{figure}[H]
\centering
\includegraphics[width=3.4 in]{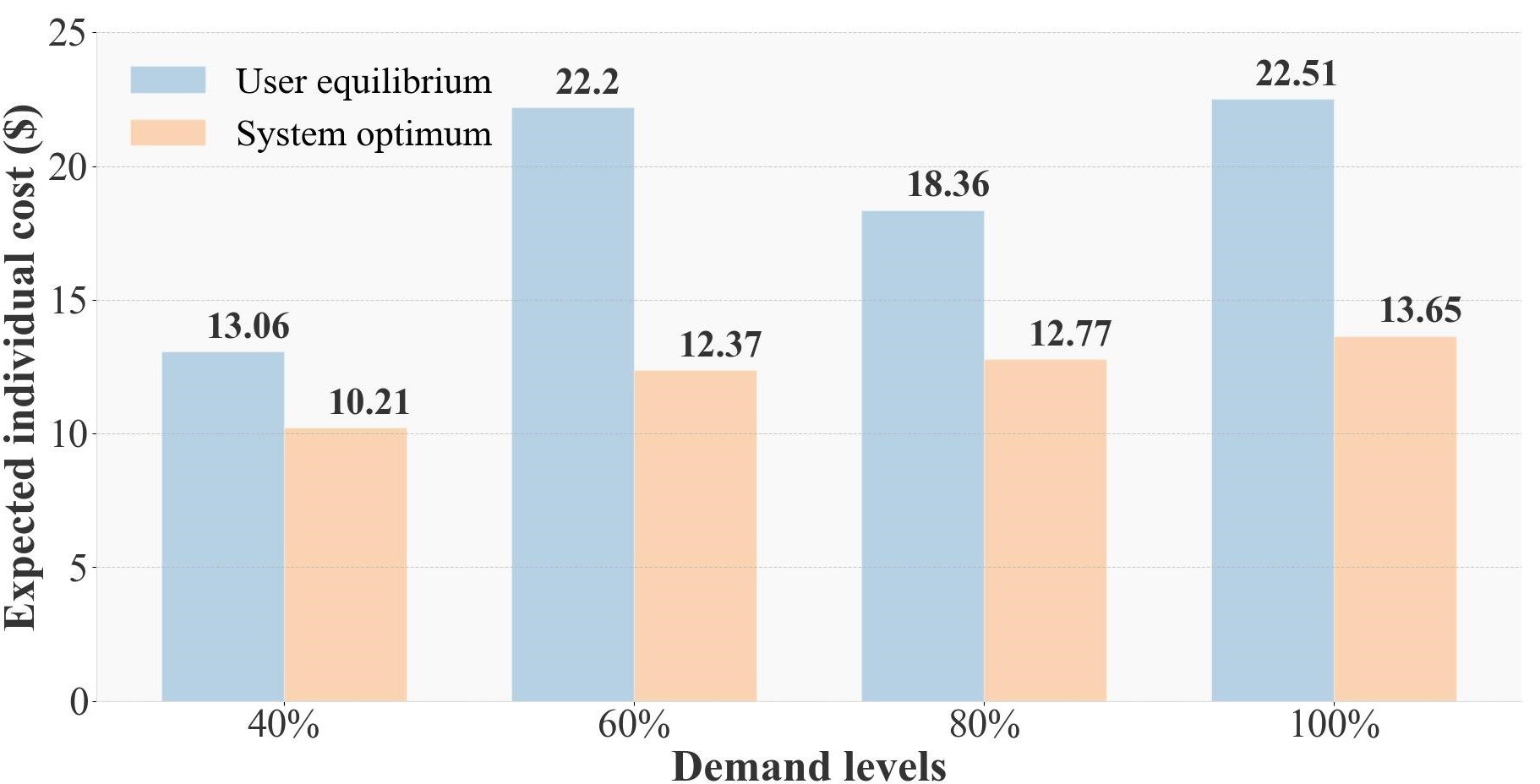}
\caption{Expected individual travel cost for OD pair `3-13' in UE and SO systems}
\label{Fig10}
\end{figure}
Taking the OD pair `3-13' as an example, Fig. \ref{Fig10} shows the corresponding expected individual travel cost (average travel cost) for the passenger. The individual cost reduction from UE to SO is ranked (from highest to lowest) as follows: scenario 60\%, scenario 100\%, scenario 80\%, and scenario 40\%. This is consistent with the potential improvement rankings: scenario 60\% (0.42), 100\% (0.36), 80\% (0.32), and 40\% (0.2). We can find the reason in Table \ref{tab:departure-arrival}, which shows that OD performance improves the most from UE to SO at the 60\% scenario. The arrival time range is optimized from (8:28 AM - 10:00 AM) to (7:28 AM - 8:28 AM), ensuring that no passenger with the OD pair `3-13' will experience a late delay, as the expected work start time is 9:00 AM. Furthermore, the travel time has also been significantly reduced from over 1 hour to less than 42 minutes.

\begin{table}[ht]
    \centering
    \caption{Departure and arrival time (AM) ranges in UE and SO systems for OD `3-13'}
    \label{tab:departure-arrival}
    \begin{tabularx}{\columnwidth}{cX X X X}
        \toprule
        \multirow{2}{*}{\textbf{Demand Level}} &  \multicolumn{2}{c} {\textbf{UE}} & \multicolumn{2}{c} {\textbf{SO}} \\
        \cmidrule(lr){2-3} \cmidrule(lr){4-5}
        & Departure & Arrival & Departure & Arrival \\
        \midrule
        100\% & 6:31-8:31 & 7:14-9:42 & 6:31-8:48 & 7:14-9:42 \\
        80\% & 6:46-8:31 & 7:28-9:42 & 6:31-8:31 & 7:14-9:27 \\
        60\% & 7:30-8:19 &  8:28-10:00 & 6:46-7:45& 7:28-8:28 \\
        40\% & 7:45-8:19 &  8:28-9:27 & 7:15-8:19 & 8:12-9:12 \\
        \bottomrule
    \end{tabularx}
\end{table}
\begin{table}[ht]
    \centering
    \caption{Interacting OD flows in the UE at demand level 60\%}
    \label{tab:intereacting-od-flow}
    \begin{tabularx}{\columnwidth}{p{0.8cm} p{1.1cm} p{1.25cm} p{1.36cm} p{2.5cm}}
        \toprule
        OD & Departure & Arrival & To station 8 & Depart after 8:00 AM \\
        \midrule
        `3-13' & 7:30-8:19 & 8:28-10:00 & 21-26 mins & 5,148 passengers\\
        `2-13' & 8:04-8:36 & 8:43-9:42 & 17 mins & 6,116 passengers\\
        `1-13' & 8:03-8:36 & 8:43-9:42 & 20 mins & 5,356 passengers \\
        \bottomrule
    \end{tabularx}
\end{table}
Focusing on OD pair `3-13', its interacting flows include all OD flows to destination 13 (`1-13', `2-13', `4-13'), as they all require a transfer at station 8, as well as flows originating from station 3 (`3-14', `3-15'). According to Fig. \ref{Fig9}, we select typical interacting OD flows from `2-13' and `1-13' to explain why the average travel cost of OD `3-13' is so high in the UE at a demand level of 60\%. Table \ref{tab:intereacting-od-flow} shows that 5,148 passengers from OD pair `3-13' (over one-third of the OD demand) will compete with passengers from OD pairs `2-13' and `1-13'. Passengers from OD pairs `2-13' and `1-13' have priority to board the train to station 8, as they are upstream flows compared to passengers from OD pair `3-13'. Therefore, more than one-third of passengers from OD pair `3-13' experience significant delays, leading to an abnormally high average travel cost for this OD.

\end{appendices}
\vfill

\end{document}